%
\input amssym.def
\input amssym.tex

\def\item#1{\vskip1.3pt\hang\textindent {\rm #1}}


\newskip\litemindent
\litemindent=0.7cm  
\def\Litem#1#2{\par\noindent\hangindent#1\litemindent
\hbox to #1\litemindent{\hfill\hbox to \litemindent
{\ninerm #2 \hfill}}\ignorespaces}

\tolerance=300
\pretolerance=200
\hfuzz=1pt
\vfuzz=1pt

\hoffset=0in
\voffset=0.5in

\hsize=5.8 true in 
\vsize=9.2 true in
\parindent=25pt
\mathsurround=1pt
\parskip=1pt plus .25pt minus .25pt
\normallineskiplimit=.99pt

\countdef\revised=100
\mathchardef\emptyset="001F 
\chardef\ss="19
\def\3{\ss}
\def\anf{$\lower1.2ex\hbox{"}$}
\def\frac#1#2{{#1 \over #2}}
\def\>{>\!\!>}
\def\<{<\!\!<}

\def\ssssarr{\hbox to 15pt{\rightarrowfill}}
\def\sssarr{\hbox to 20pt{\rightarrowfill}}
\def\ssarr{\hbox to 30pt{\rightarrowfill}}
\def\sarr{\hbox to 40pt{\rightarrowfill}}
\def\arr{\hbox to 60pt{\rightarrowfill}}
\def\larr{\hbox to 60pt{\leftarrowfill}}
\def\Arr{\hbox to 80pt{\rightarrowfill}}

\def\ad{\mathop{\rm ad}\nolimits}

\def\Aut{\mathop{\rm Aut}\nolimits}

\def\der{\mathop{\rm der}\nolimits}
\def\det{\mathop{\rm det}\nolimits}

\def\Exp{\mathop{\rm Exp}\nolimits}

\def\ev{\mathop{\rm ev}\nolimits}

\def\End{\mathop{\rm End}\nolimits}

\def\GL{\mathop{\rm GL}\nolimits}

\def\Herm{\mathop{\rm Herm}\nolimits}
\def\Hom{\mathop{\rm Hom}\nolimits}%
\def\id{\mathop{\rm id}\nolimits} 
\def\im{\mathop{\rm im}\nolimits}


%

\def\OO{\mathop{\rm O{}}\nolimits}



\def\Sym{\mathop{\rm Sym}\nolimits}

\def\Sp{\mathop{\rm Sp}\nolimits}



 
\def\UU{\mathop{\rm U{}}\nolimits}

\def\0{{\bf 0}}
\def\1{{\bf 1}}

\def\e{{\frak e}}
\def\f{{\frak f}}
\def\g{{\frak g}}
\def\gl{{\frak {gl}}}
\def\h{{\frak h}}

\def\m{{\frak m}}

\def\sp{{\frak {sp}}}

\def\uu{{\frak u}}

\def\C{{{\Bbb C}{\mskip+1mu}}} 
\def\K{{{\Bbb K}{\mskip+2mu}}} 

\def\R{{\Bbb R}} 
 
\def\N{{\Bbb N}}

\def\K{{\Bbb K}}

\def\P{{\Bbb P}} 
\def\Q{{\Bbb Q}}

\def\:{\colon}  
\def\.{{\cdot}}
\def\|{\Vert}
\def\bsk{\bigskip}

\def\giantskip{\vskip2\bigskipamount}
\def\gsk{\giantskip}
\def \la {\langle}

\def\msk{\medskip}
\def \ra {\rangle}

\def\ssk{\smallskip}

\def\bbr{\bigbreak}
\def\giantbreak{\par \ifdim\lastskip<2\bigskipamount \removelastskip
         \penalty-400 \giantskip\fi}

\def\nin{\noindent}
\def\cen{\centerline}
\def\pagebreak{\vskip 0pt plus 0.0001fil\break}
\def\linebreak{\break}

\def\eps{\varepsilon}
\def\epsilon{\varepsilon}

\def\nin{\noindent}

\def\pder#1,#2,#3 { {\partial #1 \over \partial #2}(#3)}
\def\pde#1,#2 { {\partial #1 \over \partial #2}}
\def\phi{\varphi}


\def\subeq{\subseteq}

\def\tilde{\widetilde}

\font\ninerm=cmr9
\font\eightrm=cmr8

\font\eightbf=cmbx8


\font\smc=cmcsc10
\font\bfone=cmbx10 scaled\magstep1 
\font\bftwo=cmbx10 scaled\magstep2 

\def\qed{{\unskip\nobreak\hfil\penalty50\hskip .001pt \hbox{}\nobreak\hfil
          \vrule height 1.2ex width 1.1ex depth -.1ex
           \parfillskip=0pt\finalhyphendemerits=0\medbreak}\rm}

\def\qeddis{\eqno{\vrule height 1.2ex width 1.1ex depth -.1ex} $$
                   \medbreak\rm}

\def\Lemma #1. {\bigbreak\vskip-\parskip\noindent{\bf Lemma #1.}\quad\it}

\def\Sublemma #1. {\bigbreak\vskip-\parskip\noindent{\bf Sublemma #1.}\quad\it}

\def\Proposition #1. {\bigbreak\vskip-\parskip\noindent{\bf Proposition #1.}
\quad\it}

\def\Corollary #1. {\bigbreak\vskip-\parskip\nin{\bf Corollary #1.}
\quad\it}

\def\Theorem #1. {\bigbreak\vskip-\parskip\noindent{\bf Theorem #1.}
\quad\it}

\def\Definition #1. {\rm\bigbreak\vskip-\parskip\noindent
{\bf Definition #1.}
\quad}

\def\Remark #1. {\rm\bigbreak\vskip-\parskip\noindent{\bf Remark #1.}\quad}

\def\Example #1. {\rm\bigbreak\vskip-\parskip\noindent{\bf Example #1.}\quad}
\def\Examples #1. {\rm\bigbreak\vskip-\parskip\noindent{\bf Examples #1.}\quad}

\def\Problems #1. {\bigbreak\vskip-\parskip\noindent{\bf Problems #1.}\quad}
\def\Problem #1. {\bigbreak\vskip-\parskip\noindent{\bf Problem #1.}\quad}
\def\Exercise #1. {\bigbreak\vskip-\parskip\noindent{\bf Exercise #1.}\quad}

\def\Conjecture #1. {\bigbreak\vskip-\parskip\noindent{\bf Conjecture #1.}\quad}

\def\Proof#1.{\rm\par\ifdim\lastskip<\bigskipamount\removelastskip\fi\smallskip
            \noindent {\bf Proof.}\quad}

\def\Axiom #1. {\bigbreak\vskip-\parskip\noindent{\bf Axiom #1.}\quad\it}

\def\Satz #1. {\bigbreak\vskip-\parskip\noindent{\bf Satz #1.}\quad\it}

\def\Korollar #1. {\bbr\vskip-\parskip\nin{\bf Korollar #1.} \quad\it}

\def\Folgerung #1. {\bbr\vskip-\parskip\nin{\bf Folgerung #1.} \quad\it}

\def\Folgerungen #1. {\bbr\vskip-\parskip\nin{\bf Folgerungen #1.} \quad\it}

\def\Bemerkung #1. {\rm\bigbreak\vskip-\parskip\noindent{\bf Bemerkung #1.}
\quad}

\def\Beispiel #1. {\rm\bigbreak\vskip-\parskip\noindent{\bf Beispiel #1.}\quad}
\def\Beispiele #1. {\rm\bigbreak\vskip-\parskip\noindent{\bf Beispiele #1.}\quad}
\def\Aufgabe #1. {\rm\bigbreak\vskip-\parskip\noindent{\bf Aufgabe #1.}\quad}
\def\Aufgaben #1. {\rm\bigbreak\vskip-\parskip\noindent{\bf Aufgabe #1.}\quad}

\def\Beweis#1. {\rm\par\ifdim\lastskip<\bigskipamount\removelastskip\fi
           \smallskip\noindent {\bf Beweis.}\quad}

\nopagenumbers

\def\date{\ifcase\month\or January\or February \or March\or April\or May
\or June\or July\or August\or September\or October\or November
\or December\fi\space\number\day, \number\year}

\def\title{Title ??}
\def\author{Author ??}

\def\thanks#1{\footnote*{\eightrm#1}}

\def\rightheadline{\hfil{\eightrm\title}\hfil\tenbf\folio}
\def\leftheadline{\tenbf\folio\hfil{\eightrm\author}\hfil}
\headline={\vbox{\line{\ifodd\pageno\rightheadline\else\leftheadline\fi}}}

\def\firstheadline{}
\def\firstfootline{\cen{\rm\folio}}

\def\seite #1 {\pageno #1
               \headline={\ifnum\pageno=#1 \firstheadline
               \else\ifodd\pageno\rightheadline\else\leftheadline\fi\fi}
               \footline={\ifnum\pageno=#1 \firstfootline\else{}\fi}}

\newdimen\dimenone
 \def\checkleftspace#1#2#3#4{
 \dimenone=\pagetotal
 \advance\dimenone by -\pageshrink   
 \ifdim\dimenone>\pagegoal          
   \else\dimenone=\pagetotal
        \advance\dimenone by \pagestretch
        \ifdim\dimenone<\pagegoal
          \dimenone=\pagetotal
          \advance\dimenone by#1         
          \setbox0=\vbox{#2\parskip=0pt                
                     \hyphenpenalty=10000
                     \rightskip=0pt plus 5em
                     \noindent#3 \vskip#4}    
        \advance\dimenone by\ht0
        \advance\dimenone by 3\baselineskip   
        \ifdim\dimenone>\pagegoal\vfill\eject\fi
          \else\eject\fi\fi}


\def\subheadline #1{\nin\bigbreak\vskip-\lastskip
      \checkleftspace{0.9cm}{\bf}{#1}{\medskipamount}
          \indent\vskip0.7cm\centerline{\bf #1}\medskip}
\def\subsection{\subheadline} 

\def\lsubheadline #1 #2{\bigbreak\vskip-\lastskip
      \checkleftspace{0.9cm}{\bf}{#1}{\bigskipamount}
         \vbox{\vskip0.7cm}\cen{\bf #1}\msk \cen{\bf #2}\bsk}

\def\sectionheadline #1{\bigbreak\vskip-\lastskip
      \checkleftspace{1.1cm}{\bf}{#1}{\bigskipamount}
         \vbox{\vskip1.1cm}\cen{\bfone #1}\bsk}
\def\section{\sectionheadline} 

\def\lsectionheadline #1 #2{\bigbreak\vskip-\lastskip
      \checkleftspace{1.1cm}{\bf}{#1}{\bigskipamount}
         \vbox{\vskip1.1cm}\cen{\bfone #1}\msk \cen{\bfone #2}\bsk}

\def\lchapterheadline #1 #2{\bigbreak\vskip-\lastskip\indent\vskip3cm
                       \cen{\bftwo #1} \msk \cen{\bftwo #2} \gsk}
\def\llsectionheadline #1 #2 #3{\bigbreak\vskip-\lastskip\indent\vskip1.8cm
\cen{\bfone #1} \msk \cen{\bfone #2} \msk \cen{\bfone #3} \nobreak\bsk\nobreak}


\newtoks\literat
\def\[#1 #2\par{\literat={#2\unskip.}%
\hbox{\vtop{\hsize=.15\hsize\nin [#1]\hfill}
\vtop{\hsize=.82\hsize\nin\the\literat}}\par
\vskip.3\baselineskip}

\def\references{
\sectionheadline{\bf References}
\frenchspacing

\entries\par}

\mathchardef\emptyset="001F 
\def\address{Author: \tt$\backslash$def$\backslash$address$\{$??$\}$}

\def\abstract #1{{\narrower\baselineskip=10pt{\noindent
\eightbf Abstract.\quad \eightrm #1 }
\bigskip}}

\def\addresstwo{}

\def\dlastpage{\par\vbox{\vskip1cm\nin
\line{
\vtop{\hsize=.5\hsize{\parindent=0pt\baselineskip=10pt\nin\address}}
\quad 
\vtop{\hsize=.42\hsize\nin{\parindent=0pt
\baselineskip=10pt\addresstwo}}
\hfill} }}

\def\Box #1 { \msk\par\nin 
\centerline{
\vbox{\offinterlineskip
\hrule
\hbox{\vrule\strut\hskip1ex\hfil{\smc#1}\hfill\hskip1ex}
\hrule}\vrule}\msk }

\def\adots{\mathinner{\mkern1mu\raise1pt\vbox{\kern7pt\hbox{.}}
                        \mkern2mu\raise4pt\hbox{.}
                        \mkern2mu\raise7pt\hbox{.}\mkern1mu}}


\pageno=1
\def\title{Projective completions of Jordan pairs II} 

\def\author{Wolfgang Bertram, Karl-Hermann Neeb}
\def\date{15.1.2004} 

\def\PE{\mathop{\rm PE}\nolimits}

\def\Gl{\mathop{\rm GL}\nolimits}

\def\pr{\mathop{\rm pr}\nolimits}

\def\X{\frak{X}}
\def\Aherm{\mathop{\rm Aherm}\nolimits}
\def\PP{\mathop{\Bbb P}\nolimits}

\def\address
{Wolfgang Bertram 

Institut Elie Cartan 

Facult\'e des Sciences, Universit\'e Nancy I

B.P. 239

F - 54506 Vand\oe uvre-l\`es-Nancy Cedex

France

bertram@iecn.u-nancy.fr

} 

\def\addresstwo 
{Karl-Hermann Neeb

Technische Universit\"at Darmstadt 

Schlossgartenstrasse 7

D-64289 Darmstadt 

Deutschland

neeb@mathematik.tu-darmstadt.de}

\nin
{\obeylines \parindent 0pt }
\vskip2cm
\centerline{\bfone Projective completions of Jordan pairs}
\centerline{\bfone Part II. Manifold structures and symmetric spaces} 

\gsk
\centerline{\bf\author}
\vskip1.5cm \rm

\nin {\bf Abstract.}
We define {\it symmetric spaces} in arbitrary dimension and
over arbitrary non-discrete topological fields $\K$, and we construct
manifolds and symmetric spaces associated to topological
{\it continuous quasi-inverse Jordan pairs} and {\it -triple systems}.
This class of spaces, called {\it smooth generalized projective  
geometries}, generalizes the well-known 
(finite or infinite-dimensional) bounded symmetric domains
as well as their ``compact-like'' duals.
An interpretation of such geometries as models
of Quantum Mechanics is proposed, and particular 
attention is paid to geometries that might be considered as 
``standard models"  -- they are associated
to {\it associative continuous inverse algebras} and
to {\it Jordan algebras of hermitian elements} in such an algebra.

\bigskip
\nin {\bf Contents.} 

\msk 
1. Calculus and manifolds

2. Lie groups and symmetric spaces

3. Symmetric spaces associated to continuous inverse Jordan algebras

4. Geometries associated to Jordan pairs

5. Smooth generalized projective geometries

6. Smooth polar geometries and associated symmetric spaces

7. The projective line over an associative  algebra

8. The hermitian projective line

9. Quantum mechanical interpretation

10. Prospects 

\msk \nin {\bf MSC 2000:} Prim.: 17C36, 46H70, 17C65, Sec.: 17C30, 17C90

\msk \nin {\bf Key words:} Jordan algebra, Jordan pair, Jordan triple, 
symmetric space, conformal completion, projective completion, Lie group

\sectionheadline
{Introduction}

In finite dimensions, the theory of {\it Lie groups} is closely
related to the theory of {\it symmetric spaces}.
In infinite dimensions, the theory of Lie groups is by now developed
in great generality, whereas for symmetric spaces there is not
even a commonly accepted definition. 
Nevertheless,
there is an interesting class of spaces, called {\it
(infinite-dimensional) bounded
symmetric domains}, for which one can develop a nice structure
theory and which, without doubt, are honest symmetric spaces.
Remarkably enough, the framework of their theory (developed by W.~Kaup and 
H.~Upmeier, cf.\ the monograph [Up85] and the literature given there)
is not so much Lie but rather {\it Jordan} theoretic.
Recently, also their ``compact-like" dual symmetric spaces (the analog
of the compact dual of a non-compact symmetric space in finite dimension)
have attracted attention, the most important examples being infinite-dimensional 
Gra\3mannians of many kinds (cf.\ [PS86], [DNS89], [DNS90], [KA01], [MM01], [IM02]). These
compact-like infinite-dimensional manifolds can be seen as a
``projective completion" of the underlying Jordan triple system,
in a similar way as an ordinary projective space $\R \P^n$ can
be seen as the projective completion of the affine space $\R^n$.

\ssk
In the present work, which is the second part in a series of two
papers started by [BN03], we will give a far-reaching generalization
of the above mentioned theories. We will not only free the real
theory from the Banach space set-up present in [Up85], but develop
the theory in the context of any Hausdorff topological vector space
as model space, over any non-discrete topological field. In fact,
we even work over any topological ring having dense unit group.
Compared with the approach from [Up85], our approach is more
algebraic and less analytic, which makes it considerably
simpler and more elementary. The algebraic results from Part I
of this work ([BN03]) which we need are summarized in Chapter 4,
and the basic notions of differential calculus and manifolds over
general topological fields and rings from [BGN03] are recalled in
Chapter 1. The reader who is only interested in the real or complex theory 
may everywhere replace $\K$ by $\R$ or $\C$, and he will see that all
notions from calculus we use are the ones which he is used to.

\ssk
We now give a more detailed description of the contents.
In Chapter 2 the basic theory of symmetric spaces, in arbitrary
dimension and over general base fields or rings (in which $2$ is
invertible), is developed. For several reasons, we believe
that the correct starting point for the general theory
is the approach to symmetric spaces by O. Loos ([Lo69])
-- the main idea being to incorporate all symmetries $\sigma_x$
with respects to points $x$ in the symmetric space $M$ into a
 smooth binary ``multiplication map"
$m: M \times M \to M$, $(x,y) \mapsto \sigma_x(y)$
 which is non-associative, but has 
other nice algebraic properties. The analogy with the theory of
Lie groups then becomes very close, and we get a good analog of
the functor assigning to a Lie group its Lie algebra (Theorem 2.10).
For further results on the differential geometry of 
symmetric spaces (including the {\it canonical connection} and its curvature) we
refer to [Be03b].
One should  {\it not} think of symmetric spaces as homogeneous
spaces $G/H$ -- homogeneity is a rather special phenomenon, 
and the same holds for the existence of a locally diffeomorphic {\it exponential map}
which cannot be guaranteed in general (see examples and discussion
of exponential maps in Remarks 2.11, 3.5, 6.5).

\ssk
In Chapter 3 we construct a class of symmetric spaces related
to {\it continuous inverse Jordan algebras}; by definition, these
are topological Jordan algebras over $\K$ having an open set of invertible
elements and for which the {\it Jordan inverse map} is continuous.
Once more, we closely follow the presentation from [Lo69] (cf.\ loc.\ cit.\ 
Section II.1.2.5); however, our general framework permits to
treat completely new examples such as the {\it space of non-degenerate
quadratic forms on $\K^n$} which, for fields such as $\K=\Q$,
is the prime example of a non-homogeneous symmetric space. 
For the case of Banach--Jordan algebras the symmetric space structure of the set 
of units has been studied by O.~Loos in [Lo96]. 

\ssk
Having recalled in Chapter 4 the algebraic construction and main
properties of ``generalized projective geometries" 
associated to {\it $3$-graded Lie algebras}
(which are the Lie theoretic counterpart of {\it Jordan pairs}),
we are ready to state and to prove our first main result (Theorem~5.3):
the generalized projective geometry is actually a smooth manifold
(on which the so-called {\it projective group} acts by diffeomorphisms)
if some natural conditions on the Jordan pair are fulfilled.
Namely, 
the Jordan pair $(V^+,V^-)$ shall be a topological Jordan pair over $\K$, 
the {\it set $(V^+ \times V^-)^\times$ 
of quasi-invertible pairs} shall be open in $V^+ \times V^-$, 
and the {\it Bergman-inverse mapping} 
$(V^+ \times V^-)^\times \times V^+ \times V^- \to
V^+ \times V^-$ shall be continuous; then we say that $(V^+,V^-)$
is a {\it continuous quasi-inverse Jordan pair} (Section 5.1).
If this is the case, a ``generalized quotient rule"
(Section 1.7) permits to conclude that the quasi-inverse mapping 
actually is smooth (Proposition~5.2), which is a major step in the proof
of Theorem 5.3. Our continuous quasi-inverse condition on the
Jordan pair is not only sufficient, but also necessary for
the associated generalized projective geometry to be a smooth
manifold; thus Theorem 5.3 is the most general result that one
might expect in this context. Of course, it contains the
previously mentioned results in the Banach situation as special cases. 

\ssk
In Chapter 6, we return to symmetric spaces: a symmetric space
structure on a generalized projective geometry $(X^+,X^-)$ depends on an
additional structure, namely on a fixed bijection $X^+ \to X^-$ 
which is a {\it polarity} -- in fact,
this is familiar already from the classical projective spaces
$X^+ = \R \P^n$ or $X^+=\C \P^n$: they are turned into symmetric spaces only
after the choice of a scalar product which distinguishes 
an identification of $X^+$ with the dual projective space $X^-$ and thus
determines isometry subgroups
$\P \OO_{n+1}$, resp.\ $\P \UU_{n+1}$, of the projective group 
$\P \GL_{n+1}(\K)$, $\K=\R,\C$.
We prove that, under the general assumptions of Theorem 5.3, a 
{\it continuous} polarity $p:X^+ \to X^-$ is automatically smooth and
gives rise to a symmetric space structure 
on the open set $M^{(p)}$ of non-isotropic points in $X^+$ (Theorem 6.2 (i)).
We also calculate the associated Lie triple system (i.e., the curvature
of the canonical connection; cf. [Be03b]):
it is given by anti-symmetrising the corresponding Jordan
triple product (Theorem 6.2 (ii)).
This generalizes the {\it geometric Jordan-Lie functor} which has been defined 
in [Be00] for the finite-dimensional real case.

\ssk
In Chapters 7, 8 and 9, we give applications and examples
of the preceding results and explain some links with the 
(abundant) related work in mathematics and physics.
On the one hand, Jordan algebras have been introduced by P. Jordan (cf.\
[JNW34]) in an attempt to lay algebraic foundations of quantum mechanics.
On the other hand, research on the foundations of quantum mechanics lead
by quite different arguments to the conclusion that
``... quantum mechanical systems are those whose logics form some
sort of projective geometries" ([Va85, p.\ 6]).
In the hope to bring these two lines of thought together,
the concept of ``generalized projective geometry" has been
introduced by the first named author in [Be02].
More recently, concepts of {\it delinearization of quantum
mechanics} have been proposed in the context of (Banach) hermitian
symmetric spaces, see [CGM03], where this program is motivated 
in the following way:
``The true aim of the delinearization program is to free the
mathematical foundations of quantum mechanics from any reference
to linear structure and to linear operators. It appears very
gratifying to be aware of how naturally geometric concepts describe
the more relevant aspects of ordinary quantum mechanics, suggesting
that the geometric approach could be very useful also in solving
open problems in Quantum Theories."
The close relation of the delinearization approach via hermitian
symmetric spaces to Jordan theory has not been noticed in [CGM03]
nor in the closely related paper [AS97].
In Chapter 9 we propose a ``dictionary" between
the language of generalized projective geometries (which is
equivalent to the language of Jordan theory) 
and the language of quantum mechanics.
We do not claim anything about the applicability of this 
dictionary to the ``physical world"; 
all that we aim at is to propose a terminology that makes evident 
the structural analogy between quantum mechanics and
the theory of generalized projective geometries.

\ssk
Chapters 7 and 8 are devoted to what one might call ``standard
models of quantum mechanics" -- these are the geometries corresponding
to associative continuous inverse algebras, resp.\ to their Jordan 
sub-algebras of hermitian elements. 
These are (in general)
infinite-dimensional geometries which, however, geometrically
behave very much like a projective line (over a non-commutative base
ring). A special feature of these geometries is that 
some of their associated symmetric spaces are ``of group type", i.e.
they are Lie groups, considered as symmetric spaces:
all {\it orthogonal} and {\it unitary} groups associated to involutive
continuous inverse algebras can be realized in this way.

\ssk
In the final Chapter 10 we mention some further topics and open
problems related to this work.

\msk
{\bf Notation.}
Throughout this paper, $\K$ denotes a commutative topological 
ring with unit~$1$ (i.e. $\K$ carries a topology such that the
ring operations are continuous, the group $\K^\times$
of invertible elements is open and inversion $i\: \K^\times \to \K$
is continuous) such that the group of units $\K^\times$ is
dense in $\K$.
We assume that $2$ is invertible in $\K$.
In particular, $\K$ may be any non-discrete
topological field of characteristic different from $2$
 such as $\R$, $\C$, $\Q$, $\Q_p$, $\C_p$, ${}^*\R, \ldots$

If $\K$ is a topological ring, all $\K$-modules $V$
are assumed to be topological
modules, i.e. they carry a topology such that the structure
maps $V \times V \to V$ and $\K \times V \to V$ are continuous. 
Moreover, we assume that all topological $\K$-modules are
Hausdorff. The class of continuous mappings is denoted by $C^0$.

\sectionheadline 
{1. Calculus and manifolds}

\nin {\bf 1.1. Differentiability in locally convex spaces.}
In order to motivate our general concept of differentiability,
we recall the definition of differentiable mappings on 
locally convex spaces (cf.\ [Gl01a], [Ke74], [Ha82]):
suppose $E,F$ are real locally convex spaces (not necessarily
complete), $U \subset E$
open and $f:U \to F$ continuous. Then $f$ is called 
{\it of class $C^1$} if, for all $x \in U$ and $h \in E$,
the directional derivative 
$$
df(x;h):= \lim_{t \to 0} {f(x+th)-f(x) \over t}
$$
exists and $df:U \times E \to F$ is continuous.
Inductively, one defines $f$ to be of class $C^{k+1}$ if
$df$ is of class $C^k$ (cf.\ [Gl01a, Lemma 1.14] for this
definition), and we denote by $C^0$ the class of
continuous maps.
For our purposes, the following equivalent characterization
of the class  $C^1$ will be useful:

\Proposition 1.2.
The map $f:U \to F$ is of class $C^1$ if and only if
there exists a map
$$
f^{[1]}: U \times E \times \R \supset U^{[1]}:=
\{ (x,h,t) : \, x+th \in U \} \to F
$$
of class $C^0$ such that for all $(x,h,t) \in U^{[1]}$,
$$
f(x+th)-f(x)=t \cdot f^{[1]}(x,h,t).
$$

\Proof.
Given $f^{[1]}$ as in the proposition, we get $df(x;h)=f^{[1]}(x,h,0)$,
and $df$ will be of class $C^0$ since so is $f^{[1]}$.
Conversely, assume that $f$ is $C^1$ and define $f^{[1]}$ by
$$
f^{[1]}(x,h,t):= \Big\{ \matrix{ {f(x+th)-f(x) \over t}, & \quad t \in \R^\times \cr
                         df(x)h,              & t=0. \cr}
$$
Then $f^{[1]}$ is of class $C^0$: this is seen by using,
 locally, the integral representation
$$
f^{[1]}(x,h,t)= \int_0^1 df(x+sth)h \, ds
$$
(Fundamental Theorem of Calculus, cf.\ [Gl01a, Th.\ 15]; note
that no completeness assumption is necessary here:
a priori, the integral from the right-hand side has to be taken
in the completion of $F$, but as it actually equals $f^{[1]}(x,h,t)$,
it belongs to $F$ itself.)
Now the continuity of $f^{[1]}$ follows by standard estimates 
(cf.\ [BGN03, Prop.\ 7.4] for the details).
\qed

\msk
\nin
{\bf 1.3. General definition of the class $C^1$ over
topological fields and rings.}
Now let $\K$ be a general topological ring having dense group
of units $\K^\times$, let $V,W$ be Hausdorff topological $\K$-modules
and $U \subset V$ open.
We say that a map $f:V \supset U \to W$ is $C^1(U,W)$ 
or just {\it of class $C^1$}
if there exists a $C^0$-map
$$
f^{[1]}:U \times V \times \K \supset  f^{[1]}:=
\{ (x,v,t) | \, x \in U, x+tv \in U \} \to W,
$$
such that
$$
f(x+tv)-f(x)=t \cdot  f^{[1]}(x,v,t) 
$$
whenever $(x,v,t) \in U^{[1]}$.
The {\it differential of $f$ at $x$} is defined by
$$
df(x):V \to W, \quad v \mapsto
df(x)v:=f^{[1]}(x,v,0).
$$
By density of $\K^\times$ in $\K$,
 the map $f^{[1]}$ is uniquely determined by $f$ and
hence $df(x)$ is well-defined.

\msk
\nin {\bf 1.4. Definition of the classes $C^k$ and $C^\infty$.}
Let $f\!: V \supset U \to F$ be of class $C^1$.
We say that $f$ is $C^2(U,F)$ or {\it of class $C^2$\/} 
if $f^{[1]}$ is $C^1$,
in which case we define $f^{[2]}:=(f^{[1]})^{[1]}\!:
U^{[2]}\to F$, where $U^{[2]}:=(U^{[1]})^{[1]}$.
Inductively, we say that $f$ is $C^{k+1}(U,F)$
or {\it of class $C^{k+1}$\/}
if $f$ is of class $C^k$
and $f^{[k]}\!: U^{[k]}\to F$
is of class $C^1$, in which case we
define $f^{[k+1]}:=(f^{[k]})^{[1]}\!: U^{[k+1]}\to F$
with $U^{[k+1]}:= (U^{[k]})^{[1]}$.
The map $f$ is called {\it  smooth\/} or {\it of class $C^\infty$\/}
if it is of class $C^k$ for each $k\in \N_0$.
-- Note that $U^{[k+1]}=(U^{[1]})^{[k]}$
for each $k\in \N_0$, and that
$f$ is of class $C^{k+1}$
if and only if $f$ is of class $C^1$
and $f^{[1]}$ is of class $C^k$;
in this case, $f^{[k+1]}=(f^{[1]})^{[k]}$.

\msk
\nin
{\bf 1.5. Differentiation rules.}
 We assume that $f:U \to W$ is of class $C^k$.
Its differential is the $C^0$-map
$$
df:U \times V \to W, \quad (x,v) \mapsto df(x)v=f^{[1]}(x,v,0);
$$
the {\it directional derivative in direction $v$} is
$$
\partial_v f: U \to W, \quad x \mapsto \partial_v f(x):= df(x)v.
$$
We define also
$$
Tf:U \times V \to W \times W, \quad (x,v) \mapsto (f(x),df(x)v).
$$
Then the following holds (cf.\ [BGN03]):

\ssk
\item{(1)}
For all $x \in U$,
$df(x):V \to W$ is a  $\K$-linear $C^0$-map.
\item{(2)}
If $f$ and $g$ are composable and  of class $C^k$, then $g \circ f$ 
is of class $C^k$, and $T(g \circ f)=Tg \circ Tf$.
\item{(3)}
Multilinear maps of class $C^0$ are $C^k$ and are differentiated
as usual. In particular, if $f,g:U \to \K$ are $C^1$, then the
product $f \cdot g$ is $C^1$, and
$\partial_v(fg)=(\partial_v f)g+f\partial_vg$.
Polynomial maps $\K^n \to \K^m$ are always $C^\infty$ and are
differentiated as usual.
\item{(4)}
Inversion $i:\K^\times \to \K$ is $C^\infty$, and
$(di)(x)v=-x^{-2}v$.
It follows that rational maps $\K^n \supset U \to \K^m$
are always $C^\infty$ and are differentiated as usual.
\item{(5)}
The cartesian product of two $C^k$-maps is $C^k$.
\item{(6)}
If $f:V_1 \times V_2 \supset U \to W$ is $C^1$,
and for $(x_1,x_2) \in U$ we let
$$
l_{x_1}(x_2):= r_{x_2}(x_1):= f(x_1,x_2),
$$
then the {\it rule on partial derivatives} holds:
$$
df(x_1,x_2)(v_1,v_2)= d(l_{x_1})(x_2)v_1 + d(r_{x_2})(x_1)v_2.
$$
\item{(7)} (``Schwarz' Lemma'')
If $f$ is of class $C^2$, then for all $x \in U$, $v,w \in V$,
$$
\partial_v \partial_w f(x) = \partial_w \partial_v f(x).
$$
Hence, if $f$ is of class $C^k$ and $x \in U$, then the map
$$
d^k f(x):V^k \to W, \quad
(v_1,\ldots,v_k) \mapsto \partial_{v_1} \ldots \partial_{v_k} f(x)
$$
is a symmetric multilinear $C^0$-map.
\item{(8)} There are several versions of Taylor's formula
(see [BGN03]), but none of them will be used in this work.

\msk \nin
{\bf 1.6. Continuous inverse algebras.}
We will need various generalizations of the quotient rule (4).
An associative $\K$-algebra $A$  with unit $\1$ is called a
{\it continuous inverse algebra} (c.i.a.) if the product
$A \times A \to A$ is continuous, the unit group $A^\times$ is
open in $A$ and inversion $i:A^\times \to A$ is continuous.
Writing
$$
i(x+th)-i(x)=
-x^{-1}(th)(x+th)^{-1}=
t (-x^{-1} h (x+th)^{-1}),
$$
we see that $i$ actually is $C^1$ and 
$i^{[1]}(x,h,t)=-x^{-1} h (x+th)^{-1}$, whence
$di(x)h=-x^{-1} h x^{-1}$.
Iterating this argument, we see that $i$ is $C^\infty$.

\msk
\nin {\bf 1.7.  The generalized quotient rule.}
For the second generalization of the quotient rule,
 assume $f:E \supset U \to \End(F)$
takes, on the open set $U \subset E$, values in the group
$\Gl(F)$ of (continuous) invertible linear self-maps of $E$.
We do not want to fix a topology on $\End(F)$, and hence it
makes no sense to assume $f$ or the inversion map
$j:\Gl(F) \to \Gl(F)$  to be continuous or
differentiable. Instead,
we assume that $\tilde f: U \times F \to F$,
$(x,v) \mapsto f(x)v$ is of class $C^k$ and that 
$$
\tilde{jf}:U \times F \to F, \quad (x,v) \mapsto f(x)^{-1}v
$$
is of class $C^0$. We claim that then $\tilde{jf}$ also is of class
$C^k$. Indeed, for $k=1$ we have:
$$
\eqalign{
 & 
\tilde{jf}((x,v)+s(h_1,h_2))-\tilde{jf}(x,v) \cr 
& =
\tilde{jf}((x,v)+s(h_1,h_2))-\tilde{jf}((x,v)+s(h_1,0))+
\tilde{jf}((x,v)+s(h_1,0)) - \tilde{jf}(x,v) \cr
& =
f(x+sh_1)^{-1}(v+sh_2)-f(x+sh_1)^{-1}v +
f(x+sh_1)^{-1}v - f(x)^{-1}v 
\cr
&=
s f(x+sh_1)^{-1}h_2 + (f(x+sh_1)^{-1} - f(x)^{-1}) v \cr
& =s f(x+sh_1)^{-1}h_2 +
f(x)^{-1}(f(x)-f(x+sh_1)) f(x+sh_1)^{-1} v
\cr
&=
s f(x+sh_1)^{-1}h_2+
f(x)^{-1} (\tilde f(x,f(x+sh_1)^{-1}v)-\tilde f(x+sh_1,f(x+sh_1)^{-1}v))
\cr
&=s f(x+sh_1)^{-1}h_2+
s f(x)^{-1} (\tilde f)^{[1]} ((x,f(x+sh_1)^{-1}v),(h_1,0),s)
\cr}
$$
which is the same as the product of $s$ with
$$
\eqalign{
 (\tilde{jf})^{[1]} ((x,v),(h_1,h_2),s) &
=
 f(x+sh_1)^{-1}h_2+
 f(x)^{-1} (\tilde f)^{[1]}((x,f(x+sh_1)^{-1}v),(h_1,0),s)
\cr
&=\tilde{jf}(x+sh_1,h_2)+\tilde{jf}(x, (\tilde f)^{[1]}
((x,\tilde{jf}(x+s h_1,v)),(h_1,0),s)),
\cr}
\eqno (1.1)
$$
which,
according to our assumptions, is a $C^0$-map.
It follows that $\tilde{jf}$ is $C^1$, and 
letting  $s=0$, we get
$$
d(\tilde{jf})(x,v) (h_1,h_2)=
f(x)^{-1} h_2 - f(x)^{-1} d\tilde f(x,f(x)^{-1}v) (h_1,0).
$$
Moreover, using Equation (1.1) together with the chain rule, we
can iterate this argument, and it follows that $\tilde{jf}$
is $C^k$ if so is $\tilde f$.

\msk \nin
{\bf 1.8. Manifolds.}
A {\it $C^k$-manifold with atlas (modeled on the topological $\K$-module
$E$)} 
(where $k\in \N_0\cup\{\infty\}$) is a   
topological space $M$ together with an {\it $E$-atlas\/}
${\cal A} = \{(\phi_i,U_i)\!: i \in I\}$.
This means that $U_i$, $i \in I$, is a covering of $M$ by
open sets, and
$\phi_i\!: M \supset U_i \to \phi_i(U_i)
\subset E$ is a {\it chart\/}, i.e.\
a homeomorphism of the open set $U_i \subset M$ onto an open
set $\phi_i(U_i) \subset E$, and any
two charts $(\phi_i,U_i), (\phi_j,U_j)$ are
{\it $C^k$-compatible\/} in the sense that
$$
\phi_{ij}:=
\phi_i \circ \phi_j^{-1}|_{\phi_j(U_i \cap U_j)}\! :
\phi_j(U_i \cap U_j) \to \phi_i(U_i \cap U_j)
$$
and its inverse $\phi_{ji}$ are of class $C^k$.

If the atlas $\cal A$ is {\it maximal} in the sense that it
contains all compatible charts, then $M$ is called a
{\it $C^k$-manifold (modeled on $E$)}.

{\it Smooth maps} between manifolds (with or without atlas)
are now defined as usual, and it is seen 
 that $C^k$-manifolds (with or without atlas) form
a category.

\msk \nin
{\bf 1.9. The tangent functor.}
Set-theoretically, $M$ can be seen as the quotient of the
following equivalence relation  $S/\sim$, where
$$
S:=\{ (i,x)  | \, x \in \phi_i(U_i) \} \subset I \times E,
$$
and $(i,x) \sim (j,y)$ if $\phi_i^{-1}(x)=\phi_j^{-1}(y)$.
We write $p=[i,x] \in M = S/\sim$.
Then the tangent bundle is defined to be the quotient of 
the  equivalence relation on the set
$$
TS:= S \times E\subset I \times E \times E
$$
given by:
$$
(i,x,v) \sim (j,y,w) \quad :\Longleftrightarrow \quad
\phi_j \circ \phi_i^{-1}(x)=y, \, \quad 
d(\phi_j \circ \phi_i^{-1})(x)v=w.
$$
All usual properties of the tangent bundle are now easily proved
(cf.\ [BGN03]); in particular, there is a natural manifold structure
(with atlas $T{\cal A}$) on $TM$ such that the natural
projection $\pi:TM \to M$ is smooth; the {\it tangent space} $T_p M$ is
defined to be the fiber $\pi^{-1}(p)$.
If $f\! : M \to N$ is $C^k$,
there is a well-defined  {\it tangent map\/} 
$Tf\!:  TM \to TN$, and we have the usual functorial properties
(including compatibility with direct products: $T(M \times N) \cong
TM \times TN$); thus $T$ will be called the {\it tangent functor}.

\msk \nin
{\bf 1.10. The Lie bracket.}
Smooth sections of $TM$ are called {\it vector fields}. There is a
 Lie bracket on the $\K$-module $\X(M)$ of vector fields
on $M$, given in a chart by
$$
[X,Y](x)=dY(x)X(x) - dX(x)Y(x)
\eqno (1.2)
$$
([BGN03, Th. 8.4]; note that the sign is a matter of convention). 
The Lie bracket is {\it natural} in the sense that, if 
$(X,X')$ and $(Y,Y')$ are $\phi$-related under some smooth map $\phi$,
then so is $([X,Y],[X',Y'])$ ([BGN03, Lemma 8.5]). 
See [Be03b] for a conceptual definition of the Lie bracket and for
a systematic exposition of 
differential geometry (especially, the theory of connections) 
in this framework. 

\sectionheadline
{2. Lie groups and symmetric spaces}

\nin {\bf 2.1. Manifolds with multiplication.}
A {\it product} or {\it multiplication map} 
on a manifold $M$ 
is a smooth binary map $m:M \times M \to M$, and 
{\it homomorphisms of manifolds with multiplication}
are smooth maps that are compatible with the respective 
multiplication maps.
{\it Left and right multiplication operators}, defined by
$l_x(y)=m(x,y)=r_y(x)$, are partial maps of $m$ and
hence smooth self maps of $M$. 
Applying the tangent functor to this situation, we see that
$(TM,Tm)$ is again a manifold with multiplication, and tangent
maps of homomorphisms are homomorphisms of the respective tangent
spaces.
The tangent map $Tm$ is given by the 
formula
$$
T_{(x,y)}m(\delta_x,\delta_y) =
T_{(x,y)}m((\delta_x,0_y)+(0_x,\delta_y))=
T_x (r_y) \delta_x + T_y (l_x)\delta_y.
\eqno (2.1)
$$
Formula (2.1) is nothing but the rule on partial derivatives (1.5.(6))
 written in the language of manifolds.
In particular, (2.1) shows that the canonical projection and the zero
section,
$$
\pi:TM \to M, \quad \delta_p \to p, \quad \quad
z:M \to TM, \quad p \mapsto 0_p
\eqno (2.2)
$$
are homomorphisms of manifolds with multiplication. We will always 
identify $M$ with the subspace $z(M)$ of $TM$. 
Then (2.1) implies that the operator of left multiplication by $p=0_p$
in $TM$ is nothing but $T(l_p):TM \to TM$, and similarly for
right multiplications.

\msk
\nin
{\bf 2.2. Lie groups.}
A {\it Lie group over $\K$}
is a smooth $\K$-manifold $G$ carrying a group structure such that
the multiplication map $m:G \times G \to G$ and the inversion
map $i:G \to G$ are smooth. Homomorphisms of Lie groups are smooth 
group homomorphisms. 
Clearly, Lie groups and their homomorphisms form a category in which
direct products exist.

Applying the tangent functor to the defining identities of the group
structure $(G,m,i,e)$, it is immediately seen that then $(TG,Tm,Ti,
0_{T_eG})$ is
again a Lie group such that $\pi:TG \to G$ becomes a homomorphism
of Lie groups and such that the zero section $z:G \to TG$ also
is a homomorphism of Lie groups.

\msk \nin
{\bf 2.3. The Lie algebra of a Lie group.}
A vector field $X \in \X(G)$ is called {\it left invariant} if, for
all $g \in G$,
$X \circ l_g = T l_g \circ X$. 
In particular, $X(g)=X(l_g(e))=T_e l_g X(e)$; thus $X$ is uniquely
determined by the value $X(e)$, and thus the map
$$
\X(G)^{l_G} \to T_e G, \quad X \mapsto X(e)
\eqno (2.3)
$$
from the space of left invariant vector fields into $T_e G$ is injective.
It is also surjective: if $v \in T_e G$, then
right multiplication with $v$ in $TG$, $T r_v:T G \to TG$ preserves
fibers and hence defines a vector field 
$$
v^l:G \to TG,\quad
g \mapsto T_g r_v (0_g)= Tm(g,v)=T_e l_g(v)
$$ 
which is left invariant since right
multiplications commute with left multiplications. 
Now, the space $\X(G)^{l_G}$ is a Lie subalgebra of $\X(M)$; 
this follows immediately from the naturality of the Lie bracket
because $X$ is left invariant if and only if the pair $(X,X)$ is $l_g$-related
for all $g \in G$. The space $\g:=T_eG$ with the Lie bracket defined by
$[v,w]:=[v^l, w^l]_e$ is called {\it the Lie algebra of $G$}.

\Theorem 2.4. 
\item{(i)} The Lie bracket $\g \times \g \to \g$ is $C^0$.
\item{(ii)} For every homomorphism $f:G \to H$, the tangent map
$\dot f:=T_e f:\g \to \h$ is a homomorphism of Lie algebras.

\Proof. (i)
Pick a chart $\phi:U \to V$ of $G$ such that $\phi(e)=0$. 
Since
$w^l(x)=Tm(x,w)$ depends smoothly on $(x,w)$, it is represented in
the chart by a smooth map (which again will be denoted by $w^l(x)$).
But this implies that
$[v^l,w^l](x)=d(w^l)(x) v^l(x)-d(v^l)(x)w^l(x)$
depends smoothly on $v,w$ and $x$ and hence
$[v,w]$ depends smoothly on $v,w$.

(ii)
First one has to check that the pair of vector fields 
$(v^l,(\dot \phi v)^l)$ is $f$-related, and then
the naturality of the Lie bracket implies that
$\dot f[v,w]=[\dot f v, \dot f w]$.
\qed

The functor from Lie groups over $\K$ into $C^0$-Lie algebras
over $\K$ will be called the {\it Lie functor over $\K$}.

\msk
\nin {\bf 2.5. Symmetric spaces.}
A {\it symmetric space over $\K$}
is a smooth manifold with a multiplication map $m:M \times M \to M$
such that, for all $x,y,z \in M$, writing also $\sigma_x$ for
the left multiplication~$l_x$,

\ssk
\item{(M1)} $m(x,x)=x$,
\item{(M2)} $m(x,m(x,y))=y$, i.e. $\sigma_x^2 = \id_M$, 
\item{(M3)} $m(x,m(y,z))=m(m(x,y),m(x,z))$, i.e. $\sigma_x \in \Aut(M,m)$,
\item{(M4)} $T_x(\sigma_x)=-\id_{T_xM}$.
\ssk

\nin
{\it Homomorphisms of symmetric spaces} are the corresponding
homomorphisms of manifolds with multiplication.
The left multiplication operator $\sigma_x$ 
is, by (M1)--(M3), an automorphism of 
order two fixing $x$; it is called the {\it symmetry around $x$}.
Since $2$ is invertible in $\K$, Property (M4) says that,
``infinitesimally", $x$ is an isolated fixed point of the
symmetry $\sigma_x$. If we have an implicit function theorem
at our disposition, then this holds also locally 
(see [Ne02, Lemma 3.2] for the Banach case). In particular, 
in the finite-dimensional case over $\K=\R$, 
our definition is equivalent to the one by O. Loos in [Lo69].

\ssk \nin {\bf Remark.} It would be interesting to know whether
there are real infinite-dimensional symmetric spaces for which
$x$ is not isolated in the set of fixed points of the
symmetry $\sigma_x$. 
If there were a (infinite-dimensional real) Lie group $G$ for which
the unit element is not isolated in the space of elements of order 2,
then we could take $M=G$ with $m(g,h)=gh^{-1}g$. \qed

The group $G(M)$ generated by all products $\sigma_x \sigma_y$, $x,y \in M$,
is a (normal) 
subgroup of $\Aut(M,m)$, called the {\it group of displacements}.
A distinguished point $o\in M$ is called a {\it base point}.
With respect to a base point, one defines the {\it quadratic representation}
$$
Q:M \to G(M), \quad x \mapsto \sigma_x \sigma_o.   \eqno (2.4)
$$

\Proposition 2.6.
The tangent bundle $(TM,Tm)$ of a symmetric
space is again a symmetric space.

\Proof.
We express the identities (M1)--(M3) by commutative diagrams
to which we apply the tangent functor $T$.
Since $T$ commutes with direct products, we get the same diagrams
and hence the laws (M1)--(M3) for $Tm$ (cf.\ [Lo69, II.2] for the explicit
form of the diagrams).

Next we prove (M4):
first of all, note that the fibers of $\pi:TM \to M$ (i.e. the
tangent spaces) are stable under $Tm$ because $\pi$ is a homomorphism.
We claim that for $v,w \in T_pM$ the explicit formula
$Tm(v,w)=2v-w$
holds (i.e. the structure induced on tangent spaces is the canonical
``flat'' symmetric structure of an affine space). In fact, from (M3)
for $Tm$ we get $v=Tm(v,v)=T_p(\sigma_p)v+T_p(r_p)v=-v+T_p(r_p)v$, 
whence $T_p(r_p)v=2v$ and
$$
Tm(v,w)=T_p(\sigma_p)w+T_p(r_p)v = 2v-w.
$$
Now fix $p \in M$ and $v \in T_p M$. We choose $0_p$ as base point
in $TM$.
Then $Q(v)=\sigma_v \sigma_{0_p}$ is, by (M3), an automorphism of  
$(TM,Tm)$ such that $Q(v) 0_p = \sigma_v(0_p)=2v$.
But
$$
{1 \over 2}:TM \to TM, \quad \delta_x \mapsto {1 \over 2} \delta_x
$$
also is an automorphism of $(TM,Tm)$, as shows Formula (2.1).
Therefore the automorphism group of $TM$ acts transitively on fibers,
and after conjugation of $\sigma_v$ with $({1 \over 2} Q(v))^{-1}$ we
may assume that $v=0_p$.
But in this case the proof of our claim is easy:
we have $\sigma_{0_p}=T \sigma_p$, and 
since $T_p\sigma_p = - \id_{T_pM}$, 
the canonical identification $T_{0_p}(TM) \cong T_p M \oplus T_p M$
yields
$T_{0_p}(\sigma_{0_p}) = (- \id_{T_p M}) \times (- \id_{T_pM})
= - \id_{T_{0_p} TM}$, whence (M4).
\qed

\msk \nin {\bf 2.7. The algebra of derivations of $M$.}
A vector field $X:M \to TM$ on a symmetric space $M$ is called a
{\it derivation} if $X$ is also a homomorphism of symmetric spaces.
This can be rephrased by saying that $(X \times X,X)$ is
$m$-related. The naturality of the Lie bracket therefore implies that the space
$\g$ of derivations is stable under the Lie bracket.
It is also easily checked that it is a $\K$-submodule of $\X(M)$,
and hence $\g \subset \X(M)$ is a Lie-subalgebra.

Let us fix a base point $o \in M$. 
The map $X \mapsto T \sigma_o \circ X \circ \sigma$ is a Lie algebra 
automorphism of $\X(M)$  of order 2 which stabilizes 
$\g$. We let
$$
\g=\g^+ \oplus \g^-, \quad
\g^{\pm} = \{ X \in \g \vert \, T \sigma_o \circ X \circ \sigma_o =
\pm X \}
$$
be its associated eigenspace decomposition (recall that $2$ is assumed to be
invertible in $\K$).  The space $\g^+$ is a Lie subalgebra of $\X(M)$,
whereas $\g^-$ is only closed under the triple bracket
$$
(X,Y,Z) \mapsto [X,Y,Z]:=[[X,Y],Z].
$$

\Proposition 2.8.
\item{(i)}
The space $\g^+(M)$ is the kernel of the evaluation map
$\ev_o:\g \to T_oM$, $X \mapsto X(o)$.
\item{(ii)}
Restriction of $\ev_o$ yields a bijection
$\g^- \to T_oM$, $X \mapsto X(o)$.

\Proof.
(i) Assume $X \in \g^+$.
Then $T_o \sigma X(o)=X(\sigma_o(o))=X(o)$ implies $-X(o)=X(o)$
and hence $X(o)=0$. 
On the other hand, if $X(o)=0$, then
$X(\sigma_o(p))=X(m(o,p))=Tm(X(o),X(p))=Tm(0_o,X(p))=
T\sigma_o X(p)$,
whence $X \in \g^+$.

(ii)
By (i), $\g^- \cap \ker(\ev_o) = \g^- \cap \g^+ = 0$, and hence
$\ev_o:\g^- \to T_o M$ is injective. It is also surjective:
let $v \in T_o M$. Consider the map
$$
\tilde v = {1 \over 2} Q(v) \circ z:
M \to TM, \quad p \mapsto {1 \over 2} Q(v) 0_p = {1 \over 2}
Tm(v,Tm(0_o,0_p)).
$$
It is a composition of homomorphisms and hence is itself a
homomorphism from $M$ into $TM$. Moreover, as we have seen in the proof
of Proposition 2.6, 
$\tilde v(o)=v$. Thus we will be done if we can show 
that $\tilde v \in \g^-$.
First of all, $\tilde v$ is a vector field since
$Q(v) \delta_p \in T_{m(o,m(o,p))}M=T_p M$ for all $p \in M$.
Finally,
$$
T\sigma_o \circ \tilde v \circ \sigma_o  =
{1 \over 2} T\sigma_o \circ Q(v) \circ z \circ \sigma_o 
 = {1 \over 2} Q(T \sigma_o v) \circ z =
{1 \over 2} Q(-v) \circ z = - \tilde v. 
\qeddis 

\msk
\nin {\bf 2.9. The Lie triple system of a symmetric space with
base point.}
The space $\m:=T_o M$ with triple bracket given by
$$
[u,v,w]:= - R_o(u,v)w:=[[\tilde u,\tilde v], \tilde w](o)
$$
is called the {\it Lie triple system (Lts) associated to $(M,o)$}.
It satisfies the identities of an abstract Lie triple system over
$\K$ (cf.\ [Lo69, p. 78/79]).
The notation $R_o(u,v)w$ alludes to the fact that the triple Lie bracket
indeed is the curvature tensor of a canonical connection on $M$
(cf.\ [Lo69] for the finite-dimensional real case and [Be03b]
for the general case).
Since the base point $o$ is arbitrary, we have indeed defined a 
tensor field $R$ on $M$ (in a chart it is easily seen that the
dependence of $R_o$ on $o$ is smooth).

\Theorem 2.10. Let $M$ be a symmetric space over $\K$ with base point $o$.
\item{(i)} The Lie triple bracket of the Lts $\m$ associated to $(M,o)$
is $C^0$.
\item{(ii)} If $\phi:M \to M'$ is a homomorphism of symmetric spaces 
such that $\phi(o)=o'$,
then $\dot \phi:=T_o \phi:\m \to \m'$ is an Lts homomorphism. 

\Proof.
One uses the same arguments as in the proof of Theorem 2.4.
\qed

\nin The functor from symmetric spaces with base point to 
$C^0$-Lie triple systems will be called the {\it Lie
functor for symmetric spaces}.
It contains the Lie functor for Lie groups in the following sense:
if $G$ is a Lie group, then $m(x,y)=xy^{-1}x$ defines on $G$
the structure of a symmetric space (the condition (M4) here is
equivalent to
$Ti(e)=-\id_{T_eG}$ which is proved in the same way as usual),
and as in [Lo69] it is seen that the Lts of $G$ is given in
terms of the Lie algebra of $G$ by
${1 \over 4}[[X,Y],Z]$.

\msk
\nin {\bf 2.11. On geodesics and exponential maps.} 
If $M$ is a finite-dimensional real or complex symmetric space 
and $M_1$ is a connected component of $M$, then the subgroup $G(M_1)$ of $G(M)$ 
generated by all products $\sigma_x \sigma_y$, $x,y \in M_1$, 
acts transitively on $M_1$. This follows from the existence of
an {\it exponential map} in this case (cf.\ [Lo69]).
In the general case, even for $\K=\R$, there is no exponential
map, and the connected components need no longer be homogeneous.
In the following, we give a brief account of
 the relevant definitions and explain the main arguments.

If $M$ is a symmetric space over $\K$, we define a
 {\it geodesic} to be a non-constant
homomorphism $\gamma:\K \to M$, where $\K$ carries the
``canonical flat symmetric space structure" $m(v,w)=2v-w$
which exists on any topological $\K$-module.
We say that $M$ is {\it geodesically connected}
if any two points $p,q \in M$ can be joined by a broken
geodesic, i.e. there exist points $p=p_0,\ldots,p_n=q$
such that $p_i$ and $p_{i+1}$ can be joined by a geodesic.

\Proposition 2.12. If  $M$ is geodesically connected, then
the transvection group $G(M)$ acts transitively on $M$.

\Proof. We use the same arguments as in the real finite-dimensional
case ([Lo69]):
if $\gamma:\K \to M$ is a geodesic such that
$\gamma(0)=p_i$ and $\gamma(1)=p_{i+1}$, we
let $y:=\gamma({1 \over 2})$ and
 $g:=\sigma_y \circ \sigma_{p_i} \in G(M)$; then 
$$
g(p_i)= \sigma_y(p_i)=
m(\gamma({\textstyle{1 \over 2}}),\gamma(0))=
\gamma(m({\textstyle{1 \over 2}},0))=\gamma(1)=p_{i+1}.
$$
Now the claim follows by a trivial induction on $n$.
\qed

\nin The crucial property used in the proof is that for two
points, sufficiently close to each other, we can find a
{\it midpoint}. The midpoint should be seen as a ``square root"
of one point with respect to the other; thus the lack 
of square roots in $\K$ is one obstruction for homogeneity of 
symmetric spaces, as is illustrated by the example of the projective
space $\Q \P^n$ over $\K=\Q$. Note also that geodesic connectedness does not
imply connectedness in the topological sense since already
$\K$ may be totally disconnected as the example of the
p-adic numbers $\Q_p$ shows.

We say that $M$ {\it has an exponential map} if,
for every $p \in M$ and $v \in T_p M$, there exists a unique geodesic
$\phi_v:\K \to M$ such that $\phi_v(0)=p$ and 
 $T_0 \phi_v(1)=v$ and such that the map            
$$
\Exp:=\Exp_p: T_p M \to M, \quad v \mapsto \phi_v(1)
$$
is smooth. We say that $M$ {\it is locally exponential} if $M$ has
an exponential map, and for all $p \in M$, $\Exp_p$ is a  diffeomorphism of
 some  neighborhood of $0$ in $T_pM$ onto some 
neighborhood of $p$ in $M$. Then the set of all points that can be
joined to a given point by a broken geodesic is open, and hence
$M$ is geodesically connected if $M$ is topologically connected.
It can be shown that, if $\K=\R$ and the model space of $V$
is a Banach space, then $M$ is locally exponential
(one can use the same arguments as in [Lo69]) and hence
$G(M)$ acts transitively on topological connected components.
However, already for Fr\'echet symmetric spaces this is no
longer true in general. 

\sectionheadline
{3. Symmetric spaces associated to continuous inverse Jordan algebras}

\nin {\bf 3.1. Unit groups of continuous inverse algebras.}
It is clear that the unit group $A^\times$ of a continuous inverse
algebra $A$ (cf.\ Section 1.6) is a Lie group.
The associated Lie algebra is $A$ with the commutator bracket.
We are going to explain a similar construction which arises
when one tries to replace the commutator by the anti-commutator.

\msk
\nin {\bf 3.2. Continuous inverse Jordan algebras.}
A {\it Jordan algebra} is a commutative $\K$-algebra $V$ such that the product
$x \bullet y$ satisfies the identity $x \bullet (x^2 \bullet y)=x^2\bullet
(x \bullet y)$.
Our basic reference for Jordan algebras is [MC03]; see also [FK94].
 We assume that $V$ has a unit $\1$.
Any associative algebra $A$ with the anti-commutator
$x \bullet y = {xy+yx\over 2}$ is a Jordan algebra;
subalgebras of such Jordan algebras are called {\it special}.
For $x,y$ belonging to a Jordan algebra $V$ one defines 
$$
L(x)y := x \bullet y, \quad
Q(x) := 2 L(x)^2 - L(x^2), $$
and 
$$ Q(x,y):=Q(x+y)-Q(x)-Q(y) = 2(L(x)L(y) + L(y) L(x)  - L(xy)).
\eqno (3.1)
$$
Then the {\it fundamental formula} holds:
$$
Q(Q(x)y)= Q(x) Q(y) Q(x).
\eqno (3.2)
$$
One defines the {\it Jordan inverse} $j$ by
$$
j: V^\times := \{ x \in V | \, Q(x) \,{\rm invertible} \, \} \to V, \quad
x \mapsto j(x):=x^{-1}:= Q(x)^{-1} x.
\eqno (3.3)
$$
We say that $V$ is a {\it continuous inverse Jordan algebra} (c.i.J.a.) if
$V$ is a topological Jordan algebra such that $V^\times$ is open in
$V$ and $j:V^\times \to V$ is $C^0$.

\Proposition 3.3.
The Jordan inverse of a continuous inverse Jordan algebra is smooth,
and its differential is given by
$$
dj(x)v= - Q(x)^{-1}v.
$$

\Proof. The fact that $j$ is smooth follows from the generalized
quotient rule (1.7) with $f:V \to \End(V)$, $x \mapsto Q(x)$ because
the associated map $\tilde f:(x,v) \mapsto Q(x)v$ is $C^0$ and
polynomial, hence $C^\infty$. However, in order to find the correct
expression for the differential, we repeat the main steps of the
calculation:
for $(x,h,t) \in (V^\times)^{[1]}$,
$$
\eqalign{
j(x+th)-j(x) & = Q(x+th)^{-1}(x+th) - Q(x)^{-1}x \cr
& = t Q(x+th)^{-1}h + (Q(x+th)^{-1}-Q(x)^{-1}) x \cr
&= t Q(x+th)^{-1}h -
Q(x)^{-1} (Q(x+th)-Q(x)) Q(x+th)^{-1} x \cr
&=  t Q(x+th)^{-1}h -
 Q(x)^{-1} (Q(th) + Q(x,th)) Q(x+th)^{-1} x \cr
& =
t \, 
\bigl( Q(x+th)^{-1}h - Q(x)^{-1}(t Q(h)+Q(x,h)) Q(x+th)^{-1} x \bigr). \cr}
$$
The expression following the scalar $t$ is $j^{[1]}(x,h,t)$.
Letting $t=0$, we get
$$
dj(x)h=Q(x)^{-1} h - Q(x)^{-1} Q(x,h) Q(x)^{-1} x.
$$
Now, 
$$
Q(x,h)Q(x)^{-1} x = Q(x,h)x^{-1}=
2 ([L_x,L_{x^{-1}}]+L_{xx^{-1}})h = 2 h
$$
since $L_x$ and $L_{x^{-1}}= Q(x)^{-1} L_x$ commute
(cf.\ the ``L-inverse formula'' [MC03, III.6.1]) and
$x^{-1}\bullet x=Q(x)^{-1}x^2=Q(x)^{-1} Q(x) \1 = \1$.
It follows that $dj(x)h=-Q(x)^{-1}h$.
\qed

\Theorem 3.4.
If $V$ is a continuous inverse Jordan algebra, then the set 
$M:=V^\times$ of invertible elements of $V$ is a symmetric space with product map
$$
m:M \times M \to M, \quad (x,y) \mapsto Q(x) y^{-1} =
Q(x) Q(y)^{-1} y.
$$
The quadratic map $Q:V \to \End(V)$ is a polynomial extension
of the quadratic representation $Q:M \to G(M)$ associated to the symmetric
space with base point $\1$.
The Lie triple system on the tangent space
$T_\1 M \cong V$ at the base point $\1 \in M$ is given by
$$
-R(x,y)z = [[L(x),L(y)],L(z)] \1 =  [L(x),L(y)] z = x\bullet (y \bullet z)-y
\bullet (x \bullet z).
$$

\Proof. (cf.\ [Lo96] for the Banach case) 
Using the fundamental formula (3.2), it is easily checked that
$m(x,y)$ belongs to $V^\times$ if $x,y$ belong to $V^\times$.
Thus $m$ is well-defined, and it is smooth since the Jordan
inversion is smooth (Prop.\ 3.3). 

Property (M1) follows trivially from the definition of $j$, 
(M2) and (M3) follow by straightforward calculations from
the fundamental formula (cf.\ also [Lo69, II.1.2.5]), and
since $\sigma_\1(y)=y^{-1}$, we have $\sigma_x \sigma_\1(y)=
Q(x)((y^{-1})^{-1})=Q(x)y$, proving that the quadratic
representation of $M$ and the quadratic representation of
the Jordan algebra coincide on $V^\times$.
Next we prove (M4) (using Prop. 3.3):
$$
T_x(\sigma_x) = T_x(\sigma_x \sigma_\1 \sigma_\1)=
T_x(Q(x) \circ j) = Q(x) \circ T_xj =
- Q(x)Q(x)^{-1} = - \id.
$$
In order to calculate the Lie triple system, we remark first that
$TM = T(V^\times)$ is realized by the same construction as $V$,
but with respect to the Jordan algebra
$TV \cong V \times V$ with product being the tangent map
of the Jordan product of $V$, i.e. $(x,x')\bullet (y,y')=(x \bullet y,x
\bullet y'+x' \bullet y)$
-- seen algebraically, this is the scalar extension of $V$ 
by the ring of dual numbers over $\K$, $\K[\epsilon]:=\K[x]/(x^2)
\cong \K \oplus \epsilon \K$, $\epsilon^2=0$.
Taking the unit element $\1$ as base point, the tangent vector
$v \in T_\1 M$ corresponds to the element $\1+\epsilon v \in
TV$. 
Recall from the proof of Proposition 2.8 the vector field
$$
\tilde v(p)={1 \over 2} Tm(v,Tm(0_\1,0_p)) = Tm({\textstyle{v \over 2}},
Tm(0_\1,0_p)) = 
Q({\textstyle{v \over 2}}) 0_p.
$$
With the preceding notation, $0_p = p + \eps 0$, $v=\1 + \epsilon v$, and $\tilde v$
is in the chart $V$  described by
$$
\tilde v(p)=Q(\1 + \epsilon {\textstyle{v \over 2}}).0_p=
(2 L(\1 + \epsilon {\textstyle{v \over 2}})^2 - 
L((\1 + \epsilon {\textstyle{v \over 2}})^2)) 0_p 
=(L(\1) + \epsilon L(v))p=p+\epsilon v \bullet p.
$$
In other words, in the chart $V$, $\tilde v$ is the {\it linear}
vector field given by the operator $L(v): V \to V$.
But then Formula (1.2) shows that the commutator of two linear
vector fields $L(x)$ and $L(y)$ is simply the (negative of) the usual bracket
$[L(x),L(y)]$ of endomorphisms  and hence the triple
commutator is given by $[[L(x),L(y)],L(z)]$, proving that 
$[x,y,z] =  [[L(x),L(y)],L(z)]  \1$.
{}From this the other formulas follow because $[L(x),L(y)] \1 = x \bullet y -y
\bullet x =0$.
\qed

One can prove a converse of Theorem 3.4:
a symmetric space $M$ which is open in a $\K$-module $V$ and such that
the quadratic map extends to a homogeneous quadratic polynomial,
is essentially given by the preceding construction;
in [Be00, Ch. II] (in the finite-dimensional real case) such spaces
have been called {\it quadratic prehomogeneous symmetric spaces}.

\msk
\nin {\bf 3.5. Remark on the orbit structure for the action of $G(M)$.}
In general, $M=V^\times$ is far from being homogeneous under
the action of the group $G(M)$: for instance, if $V=\Sym(n,\K)$ is the
Jordan algebra of symmetric $n \times n$-matrices over $\K$, 
then $V^\times$ is the space of non-degenerate quadratic
forms on $\K^n$. The group $G(M)$ is contained in $\Gl(V)$,
acting in the usual way on the space of forms.
It follows that the $G(M)$-orbits are included in congruence classes
of forms, and hence the orbit structure is at least as complicated
as the classification of (non-degenerate) quadratic forms over $\K$.

If $V$ is a Banach Jordan algebra over $\K \in \{ \R,\C \}$,
then $V^\times$ is a Banach symmetric space, hence is locally exponential
(Section 2.11). The exponential
map at the base point $\1$ is given by the usual exponential
series $e^v=\sum_{k=0}^\infty {v^k \over k !}$
(where the power $v^k$ is taken with respect to the Jordan product), and
the topological connected components of $V^\times$
are homogeneous under the transvection group.
It would be very interesting to understand the corresponding
situation for p-adic Banach Jordan algebras (where the 
exponential series does not converge everywhere).

\sectionheadline
{4. Geometries associated to Jordan pairs}

\nin
In this chapter we review the algebraic theory from [BN03];
results quoted without further comment can be found there.
Our basic reference for Jordan pairs is [Lo75].

\msk
\nin {\bf 4.1. Three-graded Lie algebras and Jordan pairs.}
A {\it $3$-graded Lie algebra} (over $\K$) is a Lie algebra
over $\K$ of the form $\g=\g_{1} \oplus \g_0 \oplus \g_{-1}$
such that $[\g_k,\g_l] \subset \g_{k+l}$, i.e., 
$\g_{\pm 1}$ are abelian subalgebras which are $\g_0$-modules,
in the following often denoted by $V^{\pm}$ or $\g_\pm$,
and $[\g_1, \g_{-1}] \subset \g_0$.
Then the linear map $D:\g \to \g$ with $D X = i X$ ($X \in \g_i$)
is a derivation, called the {\it grading element}, and
if $D$ is inner, $D=\ad(E)$, then the grading is called
an {\it inner $3$-grading}, and  $E$ is called an {\it Euler
operator}.
The pair $(V^+,V^-)$ together with the trilinear maps
$$
T^\pm:V^\pm \times V^\mp \times V^\pm \to V^\pm, \quad
(x,y,z) \mapsto -[[x,y],z]
\eqno (4.1)
$$
is a {\it (linear) Jordan pair} over $\K$, i.e. it satisfies the
identities, where we use the notation $T^\pm(X,Y)Z:=T^\pm(X,Y,Z)$:
$$
\eqalign{
T^\pm(X,Y,Z) & = T^\pm(Z,Y,X), \cr
T^\pm(X,Y) T^\pm(U,V,W)&=
T^\pm(T^\pm(X,Y,U),V,W) - \cr
& \quad \quad T^\pm(U,T^\mp(Y,X,V),W) +
T^\pm(U,V,T^\pm(X,Y,W)).\cr}
\eqno (4.2)
$$
Conversely, every linear Jordan pair arises in this way.

\ssk
\nin
{\bf 4.2. The projective elementary group.}
Let $(\g,D)$ be a $3$-graded Lie algebra over $\K$.
For $x \in \g_{\pm 1}$, 
the operator  $e^{\ad x} = \1 + \ad x
+{1\over 2}(\ad x)^2$ is a well defined automorphism of $\g$. 
The group generated by these operators,
$$ 
G := G(D):= \PE(\g,D) := \la e^{\ad x} \: x \in \g_{\pm 1} 
\ra \subeq \Aut(\g),
 \eqno (4.3)
$$
is called  the {\it projective elementary group of $(\g,D)$}. 
With respect to the fixed
$3$-grading, automorphisms $g$ of $\g$ will often be written
in ``matrix form"
$$
g =
 \pmatrix{g_{11} & g_{10} & g_{1,-1} \cr
g_{01} & g_{00} & g_{0,-1} \cr
g_{-1,1} & g_{-1,0} & g_{-1,-1} \cr}.
\eqno (4.4)
$$
In particular, the generators of $G$ are
represented by the following matrices (where $x \in \g_1$,
$y \in \g_{-1}$, $h \in H$):
$$
e^{\ad x}=
\pmatrix{\1 & \ad x & {1 \over 2} \ad(x)^2 \cr
0 & \1 & \ad x \cr 
0 & 0 & \1 \cr}, \quad
e^{\ad y}=\pmatrix{\1 & 0 & 0 \cr \ad y & \1 & 0 \cr 
{1 \over 2} \ad(y)^2 & \ad y& \1 \cr}
\quad 
h = \pmatrix{h_{11} & & \cr  & h_{00} & \cr  & & h_{-1,-1} \cr}.
\eqno (4.5)
$$
The subgroups $U^\pm := U^\pm(D):= e^{\ad \g_{\pm}}$ of $G$ 
are abelian and generate~$G$. 
We define the {\it automorphism group of $(\g,D)$} 
to be
$\Aut(\g,D) = \{ g \in \Aut(\g) : \, g \circ D = D \circ g \}$,
and  
we further define subgroups $H:=H(D)$ and $P^\pm:=P^\pm(D)$ of $G$ via 
$$ 
H := G(D) \cap \Aut(\g,D)
\quad \hbox{ and } \quad 
P^\pm := H U^\pm = U^\pm H. 
\eqno (4.6)
$$

\ssk
\nin {\bf 4.3. The projective completion.} From now on we assume
 that the grading derivation $D$ is inner,
$D = \ad(E)$. We denote by
$$
{\cal G}:=
\{ \ad(F) : \, F \in \g, \, \ad(F)^3 = \ad(F) \} \subset \der(\g)
\eqno (4.7)
$$
the space of all inner $3$-gradings.
By definition, the group $H$ is the stabilizer  of
$D$ in $G(D)$, and hence the homogeneous space
$M:=G(D)/H$ is just the orbit of $D$ under the action of
$G(D)$ on $\cal G$.
One shows that  $P^\pm$ is precisely the stabilizer group of the flag
$$
\f^\pm(\ad(E)) : \quad 0 \subset \f_1^\pm :=\g_{\pm 1} 
\subset \f_0^\pm:= \g_{\pm 1} \oplus \g_0 \subset \g .
\eqno (4.8)
$$
Flags of this type are called {\it inner $3$-filtrations of $\g$},
and the space of inner $3$-filtrations is denoted by $\cal F$.
The flags $o^\pm:=\f^\mp(\ad(E))$ are (for fixed $E$) called the
{\it canonical base points in $\cal F$}, and we denote by
$$
X^\pm := G(D).o^\pm \cong G(D)/P^\mp \subset {\cal F}
\eqno (4.9)
$$
their $G(D)$-orbits. The maps
$$
V^\pm \to X^\pm, \quad x \mapsto e^{\ad(x)}.o^\pm
\eqno (4.10)
$$
are injective, called the {\it projective completion of $V^\pm$}.
The reader may think of $X^\pm$ as a kind of ``manifold" modeled on
the $\K$-modules $V^\pm$: we will say that
$$
{\cal A}:= \{ (g(V^+),g): \, g \in G \}, \quad \quad
\phi_g:g(V^+) \to V^+, \quad g.x \mapsto x
\eqno (4.11) 
$$ 
is the {\it natural atlas of $X^+$}.
The chart domains $g(V^+)$ carry a natural structure of an affine
space over $\K$, depending only on the point $y:=g.o^- \in X^-$.
We then write $V_y:=g(V^+)$ and denote for $x,z \in V_y$ by
$$
\mu_r(x,y,z):=r z
\eqno (4.12)
$$
the product $rz$ in the $\K$-module $V_y$ with zero vector $x$.

\msk \nin {\bf 4.4. Transversality.}
The natural map from gradings to filtrations
${\cal G} \to {\cal F} \times {\cal F}$ and the corresponding
map $M \to X^+ \times X^- $, $gH \mapsto (gP^-,gP^+) $
are injective. Two filtrations $(\f,\e)$ are obtained from an
inner grading $\ad(E)$ if and only if they are {\it transversal}
or {\it complementary} in the sense that 
$$
\g = \f_1 \oplus \e_0, \quad
\g= \e_1 \oplus \f_0
$$
([BN03, Th.\ 3.6]);  we write then $\e \top \f$.

\msk \nin {\bf 4.5. Denominators and nominators.}
For $x \in V^+$ and $g \in \Aut(\g)$, 
we define
$$
d_g(x):=(e^{-\ad(x)} g^{-1})_{11}, \quad
c_g(x):=(g e^{\ad(x)})_{-1,-1},
\eqno (4.13)
$$
where the ``matrix coefficients" $h_{ij}$ are as in Equation (4.4).
Then
$$
d_g^+:=d_g: V^+ \to \End(V^+), \quad  \quad
c_g^+:=c_g: V^+ \to \End(V^-)
\eqno (4.14)
$$
are quadratic polynomial maps, called the {\it denominator}
and {\it co-denominator} of $g$ (w.r.t. the fixed inner grading
$\ad(E)$).
In particular, if $g=e^{\ad(w)}$, $ w \in V^-$, $x \in V^+$,
$$
\eqalign{
d_g(x) & = B_+(x,w) := \id_{V^+} + \ad(x) \ad(w) + {1 \over 4} \ad(x)^2
\ad(w)^2 \in \End(V^+)      \cr
c_g(x) & = B_-(w,x) := \id_{V^-} + \ad(w) \ad(x) + {1 \over 4}
\ad(w)^2 \ad(x)^2 \in \End(V^-) \cr}
\eqno (4.15)
$$
are called the {\it Bergman operators}. For $x \in V^+$ and 
$g \in \Aut(\g)$, we define the {\it nominator of $g$} to be 
$$
n_g(x) := \pr_1(e^{-\ad(x)} g^{-1} E)=(e^{-\ad(x)} g^{-1})_{10}.E. 
\eqno (4.16)
$$
Then $n_g:V^+ \to V^+$ is a quadratic polynomial.
In particular, for $g=e^{\ad(w)}$, $w \in V^-$,
$$
 n_g(x)=x - {1 \over 2} \ad(w)^2 x = x - Q^+(x)w.
\eqno (4.17)
$$

\Theorem 4.6. Let $g \in \Aut(\g)$ and $x \in V^+$.
Then $g.x \in V^+$ if and only if
$d_g(x)$ and $c_g(x)$ are invertible, and then
the value $g.x \in V^+$ is given by
$$
g.x = d_g(x)^{-1} n_g(x).
\qeddis 

\nin In particular, for $g=e^{\ad(w)}$, $w \in V^-$,  we get from Theorem 4.6 
$$
g(x)
=B_+(x,w)^{-1}(x - Q^+(x)w).
\eqno (4.18)
$$
In axiomatic Jordan theory, the last expression is denoted by
$x^w$ and is called the {\it quasi-inverse} (cf.\ [Lo75]).
A pair $(x,y) \in V^+
\times V^-$ is called {\it quasi-invertible} if the Bergman operators  
$B_+(x,y)$ and $B_-(y,x)$ are invertible.

\ssk
\nin {\bf 4.7. Jordan fractional quadratic maps.}
An {\it $\End(V^+)$-valued Jordan matrix coefficient
 (of type $(1,1)$, resp.\ of type $(1,0)$)}
is a map of the type
$$
q:V^\sigma \times V^\nu \to \End(V^+), \quad (x,y) \mapsto
(e^{\ad(x)}g e^{\ad(y)}h )_{11},
\eqno (4.19)
$$
where $\sigma,\nu \in \{ \pm  \}$ and $g,h \in G$, resp.\
$$
p:V^\sigma \times V^\nu \to V^+, \quad (x,y) \mapsto
(e^{\ad(x)}g e^{\ad(y)}h)_{10}E.
$$
These maps are quadratic polynomials in $x$ and in $y$,
and nominators and denominators are partial maps
of $p$ and $q$ by fixing one of the arguments to be zero. 
A {\it Jordan fractional quadratic map} is a
map of the form
$$
f:V^\sigma \times V^\nu \supset U \to V^+, \quad
(x,y) \mapsto q(x,y)^{-1} p(x,y),
$$
where $q,p$ are Jordan matrix coefficients of type (1,1),
resp.\ (1,0), and $U=\{ (x,y) \in V^\sigma \times V^\nu : \,
q(x,y) \in \Gl(V^+) \}$.

\Theorem 4.8.
The actions
$$
V^+ \times X^+ \to X^+
\quad {\it and} \quad
V^- \times X^+ \to X^+
$$
are given, with respect to all charts from the atlas $\cal A$
(cf.\ {\rm Eqn. (2.6)}), by Jordan fractional quadratic maps.
In other words, for all $g,h \in G$, the maps
$$
(v,y) \mapsto (h \circ \exp(v) \circ g).y, \quad
(w,y) \mapsto (h \circ \exp(w) \circ g).y
$$
are Jordan fractional quadratic, and the maps
$\mu_r$ are, in all charts, given by Jordan fractional quadratic
maps.
\qed

\sectionheadline
{5. Smooth generalized projective geometries}

\nin {\bf 5.1 Continuous quasi inverse Jordan pairs.}
Let $(V^+,V^-)$ be a {\it topological Jordan pair} over the
topological ring $\K$ (i.e. $V^+$, $V^-$ are topological
$\K$-modules such that the trilinear structure maps $T^+,T^-$
are $C^0$). 
If $\K=\R$ or $\C$ and 
the underlying locally convex spaces are Banach or Fr\'echet,
 then we speak of {\it Banach--}, resp., {\it Fr\'echet--Jordan pairs}.
For topological Jordan pairs  we
introduce the following two conditions:

\ssk
\item{(C1)} A topological Jordan pair is called a 
{\it continuous quasi-inverse Jordan pair}  or a {\it (C1)-Jordan pair}  
if the set  of quasi-invertible pairs,
$$
(V^+ \times V^-)^\times = \{ (x,y) \in V^+ \times V^- 
: \, B_+(x,y), B_-(y,x) \, {\rm
invertible} \, \},
$$ 
is open in $V^+ \times V^-$, and the ``Bergman inverse map" 
$$
(V^+ \times V^-)^\times \times V^+ \times V^- \to V^+ \times V^-, \quad
(x,a,v,b) \mapsto (B_+(x,a)^{-1}v,B_-(a,x)^{-1}b)
$$
is of class $C^0$. 

\item{(C2)} We say that a topological Jordan pair
 $(V^+,V^-)$ is a {\it (C2)-Jordan pair} or a {\it weak continuous
quasi-inverse Jordan pair} if, for any $a \in V^-$, the set
$$
U_a := \{ x \in V^+ : \, B_+(x,a), B_-(a,x) \, \,  {\rm invertible} \}
$$
is open in $V^+$,  and the ``partial Bergman inverse map"
$$
U_a \times V^+ \to V^+, \quad (x,v) \mapsto B_+(x,a)^{-1}v
$$
is of class $C^0$, and if the dual condition, 
with $V^+$ and $V^-$ interchanged, also holds.

\ssk
\nin
It is clear that condition (C1) implies (C2).
 For instance, Banach--Jordan pairs 
are automatically (C1) since in this case 
the operators $B(x,a)$ belong to the Banach algebra
$L(V)$ of continuous linear operators on $V$, and inversion
in the Banach algebra $L(V)$ is smooth (Banach algebras are
 special cases of  continuous inverse algebras, cf.\  1.6).

\Proposition 5.2.
In a {\rm (C1)}-Jordan pair, the quasi-inversion map
$$
(V^+ \times V^-)^\times \to V^+ \times V^-, \quad
(x,a) \mapsto (x^a,a^x) :=(e^{\ad(a)}.x,e^{\ad(x)}.a)
\eqno (5.1)
$$
is smooth, and in a {\rm (C2)}-Jordan pair, the partial maps
$$
U_a \to V^+, \quad x \mapsto e^{\ad(a)}.x, \quad \quad
U_x \to V^-, \quad a \mapsto  e^{\ad(x)}.a
$$
 are smooth.

\Proof. 
Assume $(V^+,V^-)$ satisfies Condition (C1). Following the notation
from Section 1.8, we let
$$
\eqalign{
f:V^+ \times V^- \to \End(V^+) \times \End(V^-), & \quad
(x,a) \mapsto (B_+(x,a),B_-(a,x)),         \cr
\tilde f:V^+\times V^- \times V^+ \times V^- \to
V^+ \times V^-, & \quad
((x,a),(x',a')) \mapsto f(x,a).(x',a'), \cr
jf:(V^+ \times V^-)^\times \to \Gl(V^+) \times \Gl(V^-), & \quad
(x,a) \mapsto (B_+(x,a)^{-1},B_-(a,x)^{-1}),         \cr
\tilde{j f}:(V^+\times V^-)^\times \times V^+ \times V^- \to
V^+ \times V^-, & \quad
((x,a),(x',a')) \mapsto (jf(x,a))(x',a')  \cr
 & \quad \quad \quad \quad \quad =
(B_+(x,a)^{-1}x',B_-(a,x)^{-1}a'). \cr}
$$
Then $\tilde f$ is 
a continuous polynomial, hence $C^\infty$, and by (C1),
 $\tilde{jf}$  is $C^0$.
The generalized quotient rule (Section 1.8) implies then that
 $\tilde{jf}$ is $C^\infty$.
We recall from Theorem 4.6 that for $x \in U_a$ we have 
$$ e^{\ad(a)}.x =  B_+(x,a)^{-1}(x- Q^+(x)a) \in V^+. $$ 
We therefore see that the map
$$
(x,a) \mapsto e^{\ad(a)}.x =  B_+(x,a)^{-1}(x- Q^+(x)a) 
 = \tilde{jf}(x,a,x-Q^+(x)a,0)
$$
is $C^\infty$, and that the quasi-inversion map
is $C^\infty$.
The second claim is proved by similar arguments. \qed

\Theorem 5.3. {\rm(Manifold structure on $X^\pm$)} 
Let $(V^+,V^-)$ be a topological {\rm (C2)}-Jordan pair over
the topological ring $\K$
 and $(X^+,X^-)$ its projective completion.
\item{(i)}
There exist on $X^\pm$  structures of a smooth manifolds,
modeled on the topological $\K$-modules $V^+$, resp., $V^-$,
uniquely defined by the condition  that the
collection of charts ${\cal A}^\pm=(g(V^\pm),g \in G)$
defined in {\rm Equation (4.11)} becomes an atlas of $X^\pm$.
\item{(ii)}
The projective group $G$ acts by diffeomorphisms of $X^+$ and of $X^-$.
If $g \in G$ is such that $d_g(x)$ is invertible for some
$x \in V^+$, then  the set 
$$
V_{(g)} := \{ x \in V^+ : \, d_g(x) \in \GL(V^+), c_g(x) \in \GL(V^-)
\} = \{ x \in V^+ \: g.x \in V^+\} 
$$
is open in $V^+$, and $g \: V_{(g)} \to V, x \mapsto d_g(x)^{-1}n_g(x)$ is 
a smooth map whose differential at the point $x$ is given by
$$
dg(x)=d_g(x)^{-1}.
$$
\vskip 0mm
\nin
If, in addition, $(V^+,V^-)$ satisfies {\rm (C1)}, then we have 
with respect to the manifold structure defined in Part {\rm (i)}:
\item{(iii)}
The actions 
$V^+ \times X^+ \to X^+$ and $V^- \times X^+ \to X^+$
are smooth. 
\item{(iv)}
The set $M \subset (X^+ \times X^-)$ of transversal pairs
is open in $X^+ \times X^-$.
\item{(v)}
For $r \in \K^\times$, the multiplication map
$$
\mu_r: (X^+ \times X^- \times X^+)^\top 
:= \{(x,y,z) \: (x,y), (z,x) \in M\}  \to X^+
$$
(cf.\ {\rm Equation (4.12)}) is defined on an open set and is smooth. 

\Proof.
We prove (i) for $X:=X^+$.
Uniqueness of the differentiable structure is clear since
the sets $g(V^+)$, $g \in G$, cover $X$.
In order to prove existence, we equip 
$X$ with   the final topology with respect to the maps
(the finest topology for which all these maps are continuous) 
$$
\phi_g \: V^+ \to X, \quad v \mapsto g.v,
$$
for $g \in G$, where $\phi_e$ is the inclusion
$V^+ \subset X$. In other words,
a subset $O \subeq X$ is open if and only if all
inverse images $\phi_g^{-1}(O)=g^{-1}(O) \cap V^+$,
 $g \in G$, are open in $V^+$. 

\ssk Step 1.
 $G$ acts by homeomorphisms on $X$. 
This is immediate from the definition of the topology on $X$.

\ssk Step 2.
Let us show that the induced topology on $V^+ \subset X$
is the original topology on $V^+$. Clearly, the intersection
of an open set $O$ of $X$ with $V^+$ is open in $V^+$ because 
$O \cap V^+ = \id^{-1}(O) \cap V^+$. 
Conversely, assume that $U \subset V^+$ is open in $V^+$.
We have to show that, for all $g \in G$,
 $g^{-1}(U) \cap V^+$ is open in $V^+$.
If this set is empty, we are done;
if not,  pick $x \in g^{-1}(U) \cap V^+$. 
Then $g \circ e^{\ad(x)}.0 = g.x \in U$,
and replacing $g$ by $g \circ e^{\ad(x)}$ we may assume
that $x=0$. Now, every $g \in G$ such that $g.0 \in V^+$
admits a unique decomposition
$$
g = e^{\ad(v)} h e^{\ad(w)}, \quad v \in V^+, h \in H, w \in V^-,
$$
(cf.\ [BN03, Th. 1.12 (4)]). Hence
$$
g^{-1} (U) \cap V^+ = \big(e^{-\ad(w)} h^{-1} e^{-\ad(v)}(U)\big) \cap V^+ =
\big(e^{-\ad(w)} h^{-1} (U-v)\big) \cap V^+.
$$
Now it suffices to show that $h^{-1}(U-v)$ is open in $V^+$
because $e^{\ad(w)}$, on its open domain of definition,
is smooth, hence in particular continuous 
(Proposition 5.2).  For this, we will use the following lemma:

\Lemma 5.4.
Assume $(V^+,V^-)$ is a topological {\rm (C2)}-Jordan
pair and let $B^+ \subset \End_\K(V^+)$ be the associative
subalgebra generated by all Bergman operators
$B^+(x,y)$, $x\in V^+, y \in V^-$.
Then, for all $g \in G$ and for all $x \in V^+$, the denominator 
$d_g(x)$ belongs to $B^+$.

\Proof.
We prove the lemma by induction on the ``word length of $g$''
which is, by definition, the smallest $k \in \N$ such that $g$
has an expression of the form
$$
g=e^{\ad(w_1)} e^{\ad(v_1)} \cdots e^{\ad(w_k)} e^{\ad(v_k)},
\quad v_i \in V^+, w_i \in V^-.
$$
If $k=1$, then, using the cocycle relation
$d_{fh}(x)= d_h(x) d_f(h.x)$ which holds whenever $h.x \in V^+$ (cf.
[BN03, Prop. 2.6]), we see that 
$$
d_g(x) = d_{e^{\ad(v_1)}}(x) d_{e^{\ad(w_1)}} (x+v_1)=
B(x+v_1,w_1) 
$$
belongs to $B^+$ whenever $(x+v_1,w_1)$ is quasi-invertible.
 The set of such $x$ is open in $V^+$ since our Jordan pair is (C2),
and hence generates $V^+$ as a $\K$-module. 
Therefore the denominator  $d_g:V^+ \to \End(V^+)$, being
quadratic polynomial by 4.5, coincides with the quadratic polynomial
$x \mapsto B(x+v_1,w_1)$, whence $d_g(x) \in B^+$ for all
$x \in V^+$.

Now let $g \in G$ be arbitrary  and assume
 that the claim holds for all elements of $G$ of smaller
word length than $g$.
We write $g = \tilde g \circ e^{\ad(w_k)} e^{\ad(v_k)}$
with $\tilde g$ of word length smaller than the one of $g$.
Then, again
using the cocycle relations, we have
$$
d_g(x)= d_{\tilde g e^{\ad(w_k)}}(x+v_k) =
B(x+v_k,w_k) \circ d_{\tilde g}(e^{\ad(w_k)}(x+v_k))
$$
whenever $(x+v_k,w_k)$ is quasi-invertible.
By induction, the second factor $d_{\tilde g}(e^{\ad(w_k)}(x+v_k))$
belongs to $B^+$ whenever $(x+v_k,w_k)$ is quasi-invertible.
Hence $d_g(x)$ belongs to $B^+$ whenever
$(x+v_k,w_k)$ is quasi-invertible.
As above, note that the set of such $x$
is open in $V^+$.
Thus the denominator  $d_g:V^+ \to \End(V^+)$ is a quadratic polynomial
map taking, on a non-empty open set,
values in the $\K$-module $B^+$;  hence the whole  image is  in $B^+$,
and the lemma is proved.
\qed

Note that the proof of the lemma immediately
 carries over to any Jordan pair such that each set
$U_a$, $a \in V^-$, generates $V^+$ as a $\K$-module.
However, for general Jordan pairs this property does not always
hold -- take e.g. the ring $\K[x]$, seen as a Jordan algebra
over $\K$, where the unit group is far from generating $\K[x]$
as a $\K$-module.

Now, returning to the proof of the theorem, note that elements
of $B^+$ are {\it continuous} linear operators on $V^+$ since
so are all $B(x,y)$, $x \in V^+$, $y \in V^-$. Therefore, by Lemma 5.4,
for all $h \in H$, $h_{11}=d_h(0)$ is continuous on $V^+$.
But the action of $h$ on $V^+$ is given by $h.x= h_{11} x$,
and hence $h$ acts continuously on $V^+$.
This achieves the proof of Step 2. (Note that, in particular,
we have shown that $V^+$ is open in $X$.)

\ssk Step 3. 
The transition functions are smooth. In fact,
the transition functions are
$$
\phi_{bc}=c^{-1}b:
V^+ \cap b^{-1}c(V^+) \to V^+ \cap c^{-1}b(V^+)
$$
 for $b,c \in G$. We have already seen that they are homeomorphisms.
If the intersections are non-empty, we may as above
decompose $g:=c^{-1}b$ as a product
$g=e^{\ad(v)} h e^{\ad(w)}$; the element $e^{\ad(v)}$
with $v \in V$ acts as a translation, hence smoothly,
the element $e^{\ad(w)}$ with $w \in V^-$ acts smoothly
according to Proposition 5.2, and the element $h \in H$
is a continuous linear map by Lemma 5.4 and hence
also acts smoothly. Taken together,
Step 2 and Step 3 show that $X$ is a smooth manifold.

\ssk
(ii) The proof of Step 3 above shows that elements $g \in G$
act smoothly on $X$.
 It only remains to show that the differential of $g$
is related to the denominator via  $dg(x)=d_g(x)^{-1}$. 
As above, we first reduce to the case $g \in P^-$ and $x=o$.
Then we decompose $g=h e^{\ad(a)}$, $a \in V^-$, $h \in H$.
By the chain rule and the cocycle rule for the denominators
[BN03, Th.\ 2.10], it now suffices to prove the claim for $h$ and
$\exp(a)$ separately. Since $h$ acts linearly on $V^+$ and
$d_h(o)=h^{-1}$, we are done with the first case.
As to $\exp(a)$, we have $d_{\exp(a)}(o)=B(o,a)=\id_{V^+}$.
Hence we have to show that $d \exp(a)(o)=\id$.
This follows from
$$
e^{\ad(a)}.tv - e^{\ad(a)}.o=B(tv,a)^{-1}(tv-Q(tv)a)-0=
t \, (B(tv,a)^{-1}(v-tQ(v)a)),
$$
where the last term, divided by $t$, is a $C^0$-map
of $t$ and $v$ taking value $v$ for $t=0$.

 \ssk
(iii) Recall that, according to [BN03, Th.\ 3.7], both actions are
described in charts by Jordan fractional quadratic
maps as defined in Section 4.
Therefore it suffices to show that Jordan fractional
quadratic maps are smooth: first of all, if the elements
$g,h \in G$ appearing in the definition from 4.7 are trivial,
then our claim amounts to saying that the quasi-inversion
map is smooth, which is true in a (C1)-Jordan pair,
according to Proposition 5.2. If $g$ and $h$ are not trivial,
then they can be written as a composition of translations
and quasi inverses which, according to step (i),
act as diffeomorphisms. Hence all Jordan fractional quadratic
maps are smooth.

\ssk
(iv) $M \cap (V^+ \times V^-)=(V^+ \times V^-)^\times$
is open by Property (C1).

\ssk
(v) The argument proving this claim is the same as for part (iii),
using that also $\mu_r$ is given by Jordan fractional quadratic
maps [BN03, Th.\ 4.3].
\qed

\Theorem 5.5.
Assume $(V^+,V^-)$ is a {\rm (C2)}-Jordan pair. 
Then there are canonical $G$-equivariant
bijections between the tangent bundle $TX^+$ of $X^+$ as
a smooth manifold, the tangent bundle of $X^+$ as defined in
{\rm [BN03, Th.\ 2.1]} and the tangent bundle as defined in {\rm
[Be02]} via scalar extension by dual numbers.

\Proof.
For all three models of the tangent bundle, the tangent space
$T_o X^+$, as a $\K$-module, is isomorphic to $V^+$.
Therefore in all models we get a homogeneous bundle of
the kind $G \times_{P^-} V^+$, and we only have to show
that the actions of the stabilizer group
$P^-$  on $V^+$ coincide in these three pictures.
In the context of smooth manifolds, the group $P^-$ acts
on $V^+$ via the linear isotropy representation
$\pi(p)=T_o(p)$. In the chart $V^+$, using Theorem 5.3(ii), we get 
$T_op=dp(0)=d_p(0)^{-1}=p_{11}$.
This is the representation of $P^-$ used in the model
for the tangent bundle in [BN03], and hence these
two models coincide. Finally, for the model used in [Be02],
as shown in [Be02, (7.3)],
the action of $U^-$ on $V^+$ is trivial, and
the action of $H$ commutes with $\epsilon_{o,o'}$,
so $H$ acts on $V^+$ as group of automorphisms of the
Jordan pair $(V^+,V^-)$. This characterizes the representation
of $P^-$ used in the other two models, and hence all
three models are isomorphic as $G$-bundles.       \qed

\sectionheadline{6. 
Smooth polar geometries and associated symmetric spaces}

\nin
{\bf 6.1. Continuous inverse Jordan triple systems.}
Assume $(\g,D)$ is a 3-graded Lie algebra
with an involution  $\theta$ (automorphism of order 2 reversing the grading). 
Then $V:=\g_1$ together with the trilinear map $T:V \times V \times V
\to V$ defined by
$$
T(x,y,z):= [[x,\theta(y)],z]
\eqno (6.1)
$$
is a {\it Jordan triple system} (Jts) which, by definition, is a $\K$-module
$V$ with a trilinear map $T:V \times V \times V \to V$
satisfying the identities (4.2) with superscripts omitted.
The map $V^+ \to V^-$ induced by $\theta$ is an {\it involution}
of the ``underlying Jordan pair" $(V^+,V^-)\cong (V,V)$, and in this way
Jordan triple systems are equivalent to Jordan pairs with
involution (cf.\ [Lo75]). (Note that $T$ defines a Jts if and only if $-T$
defines a Jts; thus the sign in (6.1) is a matter of convention.
Here we follow, as in [Be00], the convention that, in the real finite-dimensional
case, {\it negative} triple systems shall correspond to {\it compact} 
symmetric spaces, see below.)
A topological Jordan triple system is called (C1) or a
{\it continuous quasi inverse Jts} 
if the underlying Jordan pair $(V,V)$ is (C1) and the involution
is continuous. (For Jordan triple
systems, Condition (C2) is not very interesting.) 
Equivalently, (C1) means that
the set $(V \times V)^\times$ is open in $V \times V$
and the Bergman inverse map $(V \times V)^\times \times V \to V$
is continuous. 

\msk \nin {\bf 6.2. Polarities.}
Every involution $\theta$ of a given inner $3$-graded  
Lie algebra $(\g,D)$ induces a bijection
$$
p:X^+ \to X^-, \quad g P^- \to \theta(g) P^+.
$$
We say that $p$ is a {\it polarity} because it is an anti-automorphism 
of the generalized projective geometry $(X^+,X^-)$ (in the sense
of [Be02, Ch.\ 3]) and the corresponding
{\it space of non-isotropic points}
$$
M^{(p)} = \{ x \in X^+ : \, (x , p(x)) \in M \}
\eqno (6.2)
$$
contains the base point $o^+$ and hence is non-empty.
The {\it multiplication map}
$$
m: M^{(p)}  \times M^{(p)} \to M^{(p)}, \quad
(x,y) \mapsto \mu_{-1}(x,p(x),y)
\eqno (6.3)
$$
is well-defined and satisfies the algebraic identities
(M1)--(M3) of a symmetric space (cf.\ [Be02, 4.1], [BN03, 4.2]).
Note that, if we identify $X:=X^+$ with $X^-$ via
the polarity $p$, then by definition of $M^{(p)}$,
$$
M^{(p)} \to (X \times X)^\top \cap {\rm diag}(X \times X), \quad
x \mapsto (x,x)
\eqno (6.4)
$$
is a bijection, and hence in the chart $V=V^+ \subset X^+ = X$,
$$
M^{(p)} \cap V = \{ x \in V : \, B(x,x) \, {\rm invertible} \, \}.
\eqno (6.5)
$$
This set is open in $V$ if $(V,T)$ is (C1).

\Theorem 6.3.
Assume that $(V,T)$ is a {\rm (C1)}-Jordan triple system. 
\item{(i)}
The associated set $M^{(p)}$ of non-isotropic
points is an open  submanifold of $X$ containing the base point $o$, and 
together with the multiplication map defined by
{\rm Equation (6.3)} it is a symmetric space. 
Moreover, for all $x \in M^{(p)}$, $x$ is an isolated fixed point
of the symmetry $\sigma_x=m(x,\cdot)$.
\item{(ii)}
The Lie triple system associated to $(M^{(p)},o)$ is the vector space
$V=V^+$ together with the bracket given by
$$
[X,Y,Z]= T(X,Y,Z) - T(Y,X,Z). 
$$

\Proof.  (i)
According to Theorem 5.3 (iv), $M$ is open in
$X^+ \times X^-$. Since $p$ is $C^0$,
$M^{(p)} = \{ x \in X^+ : \, (x , p(x)) \in M \}$
is open in $X^+$.

As mentioned above,
the identities (M1), (M2), (M3) hold already in the purely
algebraic context of any generalized projective geometry (topological 
or not) with polarity. Let us prove (M4):
the involution $\sigma_x$ is given by the element $(-1)_{x,p(x)}$
of the group $G$ and hence acts as a diffeomorphism.
W.l.o.g. we may assume that $x=o$; then in the chart $V$
this diffeomorphism is given by $-\id_V$, and hence
(M4) holds. Moreover, $0$ is the only fixed point of $\sigma_o=
- \id_V$ in the open neighborhood $M \cap V$ of $o$ in $M$.

It only remains to show that $\mu$ is smooth.
This follows from the fact that $\mu(x,y)=\mu_{-1}(x,x,y)$
(when identifying $X^+$ with $X^-$),
and $\mu_{-1}$ is smooth by Theorem 5.3 (v).

(ii)
Theorem 5.5 allows us to use the realization of $TX^\pm$
from [BN03, 2.4]; in particular, we see that in the chart
$V=V^+$, vector fields $Y \in \g$ are realized by
quadratic polynomial maps $\tilde Y^+:V^+ \to V^+$.
We identify
$v \in \g_1$ with the constant vector field on $V^+$ taking value $v$.
Then $\theta(v)$ is a homogeneous quadratic vector field on $V^+$, and hence
 $\tilde v=v+\theta(v)$ is the unique vector field
in $\g^\theta$ anti-fixed by $(-\id)_*$ such that
$\tilde v(o)=v$ (here $o=o^+$). Hence the Lie triple product is given by
$$
\eqalign{
[u,v,w] & = [[\tilde u,\tilde v],\tilde w]_{o}  
= [[u+\theta(u),v+\theta(v)],w+\theta(w)]_{o} \cr
&=  [[u,\theta(v)],w]_{o} + [[\theta(u),v],w]_{o} 
 = T(u,v,w)-T(v,u,w). \cr}   
\qeddis

\msk
\nin {\bf 6.4. Remark on the orbit structure of $M^{(p)}$.} 
The space $M=(X^+ \times X^-)^\top$ is always a homogeneous
symmetric space $M \cong G/H$ with $G$ and $H$ as in 4.2,
but $M^{(p)}$, which can be seen  as the intersection of $M$
with the diagonal in $X \times X$ (cf.\ Equation (6.4)),
is in general not homogeneous under its transvection group.
A typical example for this situation is given by the
projective spaces over $\K=\Q$ or $\K=\Q_p$:
here $G/H \cong \Gl_{n+1}(\K)/\Gl_n(\K) \times \Gl_1(\K)$, but
$\K \P^n$ is not homogeneous under $\OO(n+1,\K)$.

\msk \nin
{\bf 6.5. Remark on the exponential mapping.}
Assume $V$ is a Banach Jts over $\K \in \{ \R, \C \}$.
Then $V$ is (C1), and the symmetric space $M^{(p)}$ is a Banach
symmetric space and hence is locally exponential with $\Exp=\Exp_o$.
The explicit formula for $\Exp$ is obtained as in [Be00,
Ch. X.4]: for all $x,y \in V$, the series
$$
\cosh(x)y:=\sum_{k=0}^\infty
{Q(x)^k \over (2k)! } y, \quad \quad \sinh(x):=
\sum_{k=0}^\infty {Q(x)^k x \over (2k+1) !}
$$
converge absolutely and define analytic mappings
$\cosh:V \to \End(V)$, $\sinh:V \to V$. The domain
$D:=\cosh^{-1}(\Gl(V))$ is open in $V$ and non-empty
since $\cosh(0)=\id_V$. Then, for $x \in V$, 
the exponential image $\exp(x)$ belongs to $M \cap V$
if and only if $x \in D$, and we have
$$
\exp(x)= \tanh(x):=\cosh(x)^{-1} \sinh(x) 
$$
(cf.\ [Be00, Th.\ X.4.1]; the proof carries over to the
Banach case without any changes). 
As for the case of prehomogeneous symmetric spaces (Section 3.5),
it would be very interesting to have analogous results
in the p-adic Banach case (where the series $\cosh$ and $\sinh$ do no longer
converge everywhere). 

\msk \nin
{\bf 6.6. Remark on classification.}
It goes without saying that a classification of continuous quasi
inverse Jordan pairs or -triple systems is out of reach. 
In the finite-dimensional complex or real case, {\it simple}
objects can be classified (work of O. Loos, E. Neher and others;
cf.\ [Be00, Ch. IV and XII] for precise references).
On finds that in fact essentially {\it all classical and about
half of the exceptional real and complex simple symmetric spaces 
are obtained in the form $M^{(p)}$}; this list is far too long
to be given here (see [Be00, Ch. XII]). For other base fields, so far
very little is known. In infinite dimensions over $\K\in \{ \C,\R \}$, 
various classifications
of certain simple objects are known (cf. [MC03], [Up85], [Ka83] 
(simple $JH^*$-triples) , [dlH72] (irreducible Riemannian symmetric spaces)).
--
In the following two chapters we will specialize our theory
to two important types of Jordan algebras, namely to
{\it associative algebras} and to {\it Jordan algebras of hermitian elements.}

\sectionheadline
{7. The projective line over an associative algebra}

\msk \nin {\bf 7.1. Associative algebras as Jordan pairs.}
In this chapter, $A$ is an associative algebra with unit $\1$ over
a commutative ring $\K$ having ${1 \over 2} \in \K$. 
Then $A$ is a Jordan algebra with Jordan product
$a \bullet b ={ab+ba \over 2}$ and a Jordan triple system with
triple product $T(x,y,z)=xyz+zyx$. It follows that 
the Bergman operator is given by
$$
B(x,y)z=(\1-xy)z(\1-yx)=l(\1-xy)r(\1-yx)z
\eqno (7.1)
$$
where $l(a)$ and $r(a)$ are left-, resp.\ right multiplication by $a$ in $A$.
Thus $(x,y)$ is quasi-invertible if and only if
$\1-xy$ and $\1-yx$ are invertible, and then
$B(x,y)^{-1}z=(\1-xy)^{-1} z (\1-yx)^{-1}$.
If $\K$ is a topological ring and $A$ is a continuous inverse algebra
(c.i.a.), then  the set of quasi-invertible pairs is open
in $A \times A$, and the Bergman-inverse map is continuous.
Therefore $A$ is then a (C1)-Jordan triple system and hence
$(A,A)$ is a (C1)-Jordan pair.

\msk \nin
{\bf 7.2. The three-graded picture.}
The $\K$-Lie algebra $\g:=\gl_2(A)$ of
$2 \times 2$-matrices with coefficients in $A$ has a natural
$3$-grading
$$
\g = \g_1 \oplus\g_0 \oplus \g_{-1} =
\pmatrix{0 & A \cr 0 & 0 \cr} \oplus
\pmatrix{A & 0 \cr 0 & A \cr} \oplus \pmatrix{0 & 0 \cr A & 0 \cr}
$$
which is given by the Euler operator
$$
E:={\textstyle{1 \over 2}} \, I_{1,1}:= {\textstyle{1 \over 2}}
 \pmatrix{\1 & 0 \cr 0 & - \1 \cr}.
\eqno (7.2)
$$
This $3$-grading has a natural involution given by 
$$
\theta(X)=FXF, \quad \quad F:=\pmatrix{0 & \1 \cr \1 & 0 \cr}, \quad 
\theta\pmatrix{ a & b \cr c & d \cr} = \pmatrix{ d & c \cr b & a \cr}.
\eqno (7.3)
$$
{}From the commutator relation
$$
\big[ \big[ \pmatrix{0 & x \cr 0 & 0 \cr}, \pmatrix{0 & 0 \cr y & 0 \cr} \big],
\pmatrix{0 & z \cr 0 & 0 \cr} \big] =
\pmatrix{0 & xyz + zyx \cr 0 & 0 \cr}
$$
it follows that the Jordan triple system associated to these data
is $A$ with $T(x,y,z)=xyz+zyx$. Next we are going to describe
another model of the geometry associated to this Jordan triple system.

\msk
\nin {\bf 7.3. The projective line.}
If $A$ is an associative $\K$-algebra,
we consider $W:= A \times A$ as a right $A$-module; elements
of $W$ are written as column vectors.
The {\it projective line over $A$} is, by definition, the space
$$
\P := A \P:= \{ E \subset A \times A | \, E \cong A, \exists F \cong A:
W = E \oplus F \}
$$
of $A$-submodules $E$ that are isomorphic to $A$ and admit a complement
which is also isomorphic to $A$ (cf.\ [BN03, Section 8.7], [BlHa01] or [H95]).
Elements of $\P$ can be written in the form
$$
E =\Big[ \pmatrix{x \cr y \cr} \Big] := \Big\{ \pmatrix{xa \cr ya \cr} | \, a \in A \Big\}
$$
where $\bigl( {x \atop y} \bigr)$ is a base vector of $E$ over $A$.
For $(E,F) \in \PP \times \PP$ we write $E \top F$ if
$W = E \oplus F$, and we let
$(\P \times \P)^\top = \{ (E,F) \in \P \times \P | \, E \top F \}$.
Then the map
$$
{\cal P}:= \{ p \in \End_A(W) | \, p^2 = p, \im(p) \cong A, \ker(p) \cong A \}
 \to (\P \times \P)^\top, \quad
p \mapsto (\ker(p),\im(p))
$$
is a bijection. As ``canonical'' base point in $(\P \times \P)^\top$ we choose
$(o^+,o^-)=(A \times 0,0 \times A)=
([\bigl({\1 \atop 0}\bigr)],[\bigl({0 \atop \1}\bigr)])$
 which corresponds to the projection
$p=\bigl({\1 \atop 0}{0 \atop 0}\bigr)$.
The group $\Gl_2(A)$ acts transitively on the projective line $\P$
and on the set $(\P \times \P)^\top$.
Another base point  is given by
$[\bigl({\1 \atop \1}\bigr)],[\bigl({\1 \atop -\1}\bigr)]$.
 The matrix transforming the canonical base point
into the new one is the {\it Cayley transform}
$$
C= \pmatrix{\1 &  \1 \cr \1 & -\1 \cr}.
\eqno (7.4)
$$

\msk \nin {\bf 7.4. Affine charts of the projective line.}
Every pair $(E,F) \in (\P \times \P)^\top$
defines a linearization of $\P$: the set $F^\top$ of elements that
are transversal to $F$ is an affine space over $\K$ (not over $A$
in general), and taking $E$ as origin,  $F^\top$ is turned into
a $\K$-module. This module is (non-canonically) isomorphic to $A$.
For the canonical base point $(o^+,o^-)$ we fix such an imbedding of
$A$ into $\P$:
$$
\Gamma: A \to \P, \quad z \mapsto \Gamma_z:=\Big[ \pmatrix{z \cr \1 \cr}\Big].
$$
Note that $\Gamma_z$ is the graph of the left translation
$l_z:A \to A$, $a \mapsto za$. In this picture, the action of 
$\Gl_2(A)$ is described by usual fractional linear transformations,
$$
\pmatrix{a & b \cr c & d \cr} \Gamma_z =
\Gamma_{(az+b)(cz+d)^{-1}}
\eqno (7.5)
$$
if $cz+d$ is invertible. 
In particular, the matrix $F$ from Equation (7.3) represents inversion
in $A$, and $I_{1,1}$ (Equation (7.2))
represents multiplication by the scalar $-1$.
The imbedding $\Gamma:A \to \P$ does not only depend on
the base point $(o^+,o^-)$ but also on the fixed normalization
of its representatives; however, the {\it sets} $\Gamma(A)$ and
$\Gamma(A^\times)$ depend only on $(o^+,o^-)$.
For $\Gamma(A^\times)$ a more intrinsic description is given by
$$
\Gamma(A^\times) = \Gamma(A) \cap \{ E \in \P | \, E \top I_{1,1}(E) \},
\eqno (7.6)
$$
and the projective transformation induced by $I_{1,1}$ indeed
depends only on $(o^+,o^-)$ (in fact, we have seen above that $I_{1,1}$
is induced by multiplication by the scalar $-1$ in the $\K$-module
defined by the pair $(o^+,o^-)$ and hence its effect on $\P$ depends
only on $(o^+,o^-)$). 
Moreover, the symmetric space structure on $\Gamma(A^\times)$
also depends only on the pair $(o^+,o^-)$, whereas the group
structure cannot be defined in terms of $(o^+,o^-)$ alone.

\msk \nin
{\bf 7.5. Imbedding of the projective line into the three-graded model.}
For every projection $p:W \to W$, $\ad(p):\g \to \g$
is an inner 3-grading, and for every $E = \im(p) \in \P$, we get
the corresponding flag $(\f_1 \subset \f_0 \subset \g) \in {\cal F}$
which only depends on $E$. This defines a commutative diagram of maps
$$
\matrix{ {\cal P} \cong (\P \times \P)^\top & \subset & \P \times \P \cr 
\downarrow & & \downarrow \cr
{\cal G} \cong ({\cal F} \times {\cal F})^\top & \subset & {\cal F}
\times {\cal F} \cr}
\eqno (7.7)
$$
which are all $\Gl_2(A)$-equivariant, and the vertical arrows are
injective ([BN03, Theorem 8.4]). In particular, the natural map
$\P \to {\cal F}$ is an injection, and it is a bijection when
restricted to the ``(geometric) connected components of the base
point" which are the orbits of the respective base points
under the elementary projective group $G=\P E_2(A)$, where
$$
\eqalign{
E_2(A) & =  \la P^+,P^- \ra  \quad \subset \Gl_2(A), \cr
P^+ & = \Big\{ \pmatrix{\1 & x \cr 0 & \1 \cr}| \, x \in A \Big\}, \quad
P^- = \Big\{ \pmatrix{\1 & 0 \cr y & \1 \cr}| \, y \in A \Big\}. \cr}
\eqno (7.8)
$$
Note that the matrix
$$
J:= \pmatrix{0 & \1 \cr -\1 & 0 \cr} = \pmatrix{\1 & \1 \cr 0 & \1 \cr}
\pmatrix{\1 & 0 \cr - \1 & \1 \cr} \pmatrix{\1 & \1 \cr 0 & \1 \cr}
\eqno (7.9)
$$
belongs to $E_2(A)$ and satisfies $J.o^+ = o^-$. It follows that in both 
models we have $X^+=X^-$ as sets. Moreover, since all base points in
$(\P \times \P)^\top$ are conjugate under $\Gl_2(A)$, the same results
hold also for all other geometric connected components of $\P$.

\msk \nin {\bf 7.6. Manifold structures.}
Now assume that $\K$ is a topological ring and
$A$ is a c.i.a. over $\K$. 
As we have seen in Section 7.1, $A$ is then a (C1)-Jordan triple
system, and hence the projective completion $X^+ \cong G/P^-$
of $A$ carries a natural manifold structure satisfying all
properties from Theorem 5.3.  Using the imbedding from Section 7.5,
by transport of structure, the component
$G.o^- \subset \P$ can be equipped with the same structure, and since
$\P$ is a disjoint union of geometric connected components which are
conjugate under $\Gl_2(A)$, we get a natural manifold structure on all  
of $\P$. This manifold structure agrees with the one that is obtained
by taking $\Gamma(A) \subset \P$ as ``base chart" and then 
constructing directly, via the action of $\Gl_2(A)$, 
an atlas on $\P$ in the same way as we did for $X^+$ in Chapter~5. 
This is an immediate consequence of the $\GL_2(A)$-equivariance of the 
diagram (7.7). 

\msk \nin {\bf 7.7. Symmetric space structures.}
Associated to the given base point $(o^+,o^-) \in \P \times \P$,
there are three natural involutions of $G$, given by conjugation
with the matrices $I_{1,1},F,J$, respectively.
The first two are related to each other via the Cayley transform
$C$ and give rise to the
symmeric space $A^\times \cong \Gamma(A^\times) \subset \P$.
The third one gives rise the ``c-dual symmetric space of $A^\times$"
which is isomorphic to $A[i]^\times / A^\times$, where
$A[i] := A\otimes_\K (\K[x]/(x^2 + 1))$ is the ``complexification"
of $A$ (cf.\ [Be00]).

\sectionheadline
{8. The hermitian projective line}

\msk \nin
{\bf 8.1. The space of hermitian elements.}
Assume $A$ is as in Section 7.1 and $*:A \to A$ is an involution
($\K$-linear antiautomorphism of order $2$). We define the spaces of
{\it hermitian}, resp.\ {\it anti-hermitian elements} by
$$
\Herm(A,*):= \{ a \in A | \, a^* = a \}, \quad \quad \quad
 \Aherm(A,*):= \{ a \in A | \, a^* = -a \}.
$$
Then $\Herm(A,*)$  is a Jordan-subalgebra of $A$, and $\Aherm(A,*)$
is a Jordan-sub triple system of $A$.
Recall that $2$ is assumed to be invertible in $\K$, so
$A = \Herm(A) \oplus \Aherm(A)$.
(If $\K=\R$ and $A$ is an algebra over  $\C$ such
that $*$ is $\C$-anti-linear, then $i \Aherm(A,*) = \Herm(A,*)$;
more generally, this holds whenever there is an element $j \in Z(A)$
such that $j^2 = \pm \1$ and $j^*=-j$.)
 We are going to describe Linear
Algebra models for the geometries associated to the Jordan
pairs $(\Herm(A,*),\Herm(A,*))$ and $(\Aherm(A,*),\Aherm(A,*))$. 
They will be closely related to the {\it $*$-unitary group}
$$
\UU(A,*):= \{ a \in A^\times| \, a^{-1} = a^* \}.
$$

\msk \nin
{\bf 8.2. The $*$-symplectic and the $*$-pseudo unitary group.}
If $*$ is an involution of $A$, then
by a direct calculation one checks that the $\K$-linear map
$$
\Phi_1: M_2(A) \to M_2(A), \quad
\pmatrix{a & b \cr c & d \cr} \mapsto
\pmatrix{d^* & - b^* \cr - c^* & a^* \cr}
\eqno (8.1)
$$
is an involutive anti-automorphism of the associative algebra
$M_2(A)$. If $A$ is commutative and $*=\id$, then $\Phi_1(A)$ is the 
matrix $\tilde A$ adjoint to $A$ via the relation $A \tilde A =
\tilde A A = \det(A) \1$, and then the map $\Phi_1$ appears also as 
``symplectic involution" in
the context of the Cayley--Dickson process, cf.\ [MC03, II.2.9].
We can define three other involutions of $M_2(A)$ by
$$
 \Phi_2(X):=I_{1,1} \Phi_1(X) I_{1,1}, \quad  \quad
\Phi_3(X):=F \Phi_1(X) F, \quad \quad
\Phi_4(X):=J \Phi_1(X) J^{-1}. 
\eqno (8.2)
$$
With $X=\bigl({a \atop c}{b\atop d}\bigr)$, the explicit formulas are:
$$
\Phi_2(X) =\pmatrix{d^* & b^* \cr c^* & a^* \cr}, \quad
\Phi_3(X) =\pmatrix{a^* & -c^* \cr -b^* & d^* \cr}, \quad
\Phi_4(X) =\pmatrix{a^* & c^* \cr b^* & d^* \cr}. \quad
\eqno (8.3)
$$
If $\Phi = \Phi_j$, $j=1,2,3,4$,  is any of these involutions, we
obtain an involutive automorphism of $\Gl_2(A)$ by
$$
\tilde \Phi_j: \Gl_2(A) \to \GL_2(A), \quad g \mapsto \Phi_j(g)^{-1}
$$
and an involutive Lie algebra automorphism
$$
\dot \Phi_j: \gl_2(A) \to \gl_2(A), \quad X \mapsto - \Phi_j(X).
$$
We define the {\it $*$-symplectic} and the {\it $*$-pseudo unitary
group} via
$$
\eqalign{
\Sp(A,*) & := \UU(A \times A, \Phi_1) = \{ g \in\GL_2(A) : \, \Phi_1(g)=g^{-1} \}, 
\cr
\UU(A, A, *) &
:=\UU(A \times A, \Phi_2)= \{ g \in \GL_2(A) : \,\Phi_2(g)=g^{-1} \}, \cr}
$$
and the corresponding Lie algebras 
$$
\eqalign{
\sp(A,*) & := \{ X \in\gl_2(A) : \, \Phi_1(X)= - X \}, 
\cr
\uu(A, A, *) &
:=\{ X \in \gl_2(A) : \,\Phi_2(X)=-X \}. \cr}
$$
Since $-\Phi_j(I_{1,1})=I_{1,1}$ for $j=1,2$, these two Lie algebras
are stable under the grading derivation $\ad(I_{1,1})$ of $\gl_2(A)$
and hence are themselves 3-graded Lie algebras which, moreover, are
stable under conjugation by the matrix $F$. It follows that the
Jordan triple system corresponding to these involutive $3$-graded
Lie algebras is given by restricting the one from $\gl_2(A)$ to
$-\Phi_j$-invariants.  Now,
$$
- \Phi_1 \pmatrix{0 & x \cr 0 & 0 \cr}=\pmatrix{0 & x^* \cr 0 & 0 \cr}, \quad
\quad - \Phi_2 \pmatrix{0 & x \cr 0 & 0 \cr}=\pmatrix{0 & -x^* \cr 0 & 0 \cr}, 
$$
and hence the Jts associated to $\sp(A,*)$ is $\Herm(A,*)$ and
the Jts associated to $\uu(A,A,*)$ is $\Aherm(A,*)$.

\msk \nin {\bf 8.3. The (anti-) hermitian projective line.}
Next we are going to describe the geometries associated to
$\Herm(A,*)$ and to $\Aherm(A,*)$. We have to extend the involutions $*$
and $-*$ of $A$ to globally defined maps $\P \to \P$. The idea is
simply to send an element $\im(p) \in \P$, where $p \in {\cal P}$, to
the element $\ker(\Phi_j(p))$, $j=1,2$. This is well-defined:

\Lemma 8.4. Let $V$ be a right $A$-module and
$\cal R$ be the set of all complemented  right $A$-submodules of $V$
and assume $\phi:\End_A(V) \to \End_A(V)$ is a $\K$-linear anti-automorphism.
Then the map
$$
\tilde \phi: {\cal R} \to {\cal R}, \quad \im(p) \mapsto \ker(\phi(p))
$$
(where $p$ is a projection onto $\im(p)$)
is a well-defined bijection satisfying 
$$ 
\tilde\phi(g.E) = \phi(g)^{-1}.\tilde\phi(E), \quad g \in \GL_A(V) 
= \End_A(V)^\times.
 $$
Moreover, if $V=A \times A$ and
$\phi=\Phi_j$, $j=1,2,3$, then $\P$ is stable under $\tilde \phi$.

\Proof.
First of all, if $p^2=p$, then also $(\phi(p))^2 = \phi(p)$,
hence $\phi(p)$ is a projection.
If $p$ and $q$ are projections such that $\im(p)=\im(q)$,
then there exists $g \in \Gl_A(V)$ such that $q=p \circ g$.
Hence $\ker(\phi(q)) = \ker(\phi(g) \circ \phi(p))=
\ker(\phi(p))$ since $\phi(g)$ is bijective. 
Thus $\tilde \phi$ is well-defined.
Clearly, $\tilde \phi$ is bijective with inverse $\tilde{\phi^{-1}}$.

The transformation property under $g$ follows from
$$
\phi(g)^{-1}.\ker (\phi(p)) = \ker (\phi(g)^{-1} \phi(p) \phi(g)) = 
\ker(\phi(gpg^{-1})) = \tilde \phi (\im (gpg^{-1}))= \tilde \phi(g.\im(p)).
$$

Now let $\phi=\Phi_j$, $j=1,2$, and $\im(p) \in \P$.
Then there exists $g \in \GL_2(A)$ with $gpg^{-1}=
\bigl({\1 \atop 0}{0\atop 0}\bigr)$, whence
$\phi(p)=\phi(g) \phi\bigl({\1 \atop 0}{0\atop 0}\bigr) \phi(g)^{-1} 
=\phi(g) \bigl({0 \atop 0}{0\atop \1}\bigr) \phi(g)^{-1}$,
which has kernel $\phi(g)(A \times 0) \cong A$. 
For $j=3$, it suffices to note that the matrices $F$ and $I_{1,1}$ are
conjugate to each other (cf.\ Section 8.5), and hence also
$\tilde \Phi_2$ and $\tilde \Phi_3$ are conjugate to each other.
\qed

\msk \nin
The Lemma shows  that $\tilde \Phi_j$ for $j=1,2$
 is induced by
the automorphism $\tilde \Phi_j:E_2(A) \to E_2(A)$ (which is well-defined
since, by  (8.1),  the unipotent groups $P^\pm$ 
defined in (7.8) are stable under $\tilde \Phi_j$, $j=1,2$), i.e.\ $\tilde \Phi_j$
is given by
$$
\tilde \Phi_j: \P \to \P,  \quad g.o^+ \to \tilde \Phi_j(g).o^+. 
\eqno (8.4)
$$
We say that an element $E \in \P$ is
\ssk
\item{--} {\it hermitian} if $\tilde \Phi_1(E)=E$,
\item{--} {\it anti-hermitian} if $\tilde \Phi_2(E)=E$,
\item{--} {\it unitary} if $\tilde \Phi_3(E) = E$.

\ssk \nin
Assume $E=\Gamma_z=[ \bigl( {z \atop \1} \bigr) ] = \im(p)$ with
$p= \bigl({0 \atop 0} {z \atop \1} \bigr)$.
Then $\Phi_1(p)=\bigl({\1 \atop 0} {-z^* \atop 0} \bigr)$
has kernel $[\bigl( {z^* \atop \1}\bigr)]$.
Thus the restriction of $\tilde \Phi_1$ to $A=\Gamma(A)$ is the involution $*$,
and $E$ is hermitian if and only if $z$ is hermitian.
Similarly, we see that $E$ is anti-hermitian if and only if
$z$ is anti-hermitian.
Finally, $E$ is unitary if and only if
$[\bigl( {z \atop \1} \bigr) ] = [\bigl( {\1 \atop z^*} \bigr) ]$.
First of all, this implies that $z$ must be invertible in $A$,
and then the condition
$[\bigl( {z \atop \1} \bigr) ] = [\bigl( {(z^*)^{-1} \atop \1} \bigr) ]$
is equivalent to $z^{-1} = z^*$, i.e. to the unitarity of $z$.

\ssk 
The sets $\P_h$, resp.\ $\P_{ah}$ of hermitian, resp., (anti-)hermitian 
elements in $\P$ is called
the {\it (anti-)hermitian projective line}; the set $\P_u$ of 
unitary elements is called the {\it unitary projective line}.
The projective completion of $\Herm(A,*)$, resp.\ of $\Aherm(A,*)$ 
are the imbeddings
$$
\Gamma: \Herm(A,*) \to \P_h, \quad \quad \Gamma:\Aherm(A,*) \to \P_{ah}.
$$
This geometric picture can be imbedded into the three-graded picture
simply by restricting the imbedding (7.7) to $\tilde \Phi_j$-invariants.

\msk \nin {\bf 8.5. The Cayley transform.}
The matrices $F$ and $I_{1,1}$ are conjugate in $\GL_2(A)$ via $C$:
$F=C^{-1} I_{1,1} C$.  
It follows that $C^{-1}(\P_{ah})=\P_u$, 
i.e. the anti-hermitian and the unitary projective line are
isomorphic. In particular, the unitary group $U(A,*)$ is
injected into $\P_{ah}$ via
$$
U(A,*) \to \P_{ah}, \quad z \mapsto C(\Gamma_z)  
$$
If $z-e$ is invertible, then the last term equals
$\Gamma_{(z+e)(z-e)^{-1}}$ and it belongs to $\Gamma_{\Aherm(A)}$.

\msk
\nin {\bf 8.6. Manifold structures and symmetric spaces.}
If $A$ is a c.i.a.\  over a topological ring $\K$ and $*$ is continuous,
then $\Herm(A,*)$ and $\Aherm(A,*)$ are (C1)-Jordan triple systems.
The corresponding manifold structure on the geometric models is
again simply obtained by seeing everything as submanifolds fixed
under $\tilde \Phi_j$ in the models corresponding to $A$.
The natural polarities  given by the matrix $I_{1,1}$, resp.\ by $F$,
define symmetric spaces: as explained in the preceding section,
the unitary group arises as the space of non-isotropic points in
the anti-hermitian projective line $\P_{ah}$; in particular, 
$\UU(A,*)$ is a symmetric space. Moreover, with respect to the
underlying manifold structure, also the group multiplication in
$\UU(A,*)$ is smooth (the calculation is exactly the same as the
one for the orthogonal group $\OO_n(\R)$ in the Cayley chart),
 and hence $\UU(A,*)$ is a Lie group.
The natural symmetric space realized
 in the hermitian projective line is the space of invertible
elements in the Jordan algebra $\Herm(A,*)$ (already encountered in
Chapter 3), resp.\ its c-dual symmetric space.
Since the set-up is almost the same as the one in [Be96] (where the 
special case $A=\End(V)$, $*=$ adjoint, was considered; cf.\ Example 8.7),
we can  refer to [Be96] and
to [Be00] for further details of the calculations.

In general, there are many other polarities which are not
isomorphic to the natural ones, and hence there are other symmetric
spaces that can be realized inside $\P_h$ or $\P_{ah}$. In [Be00,
XI.5] they have been called {\it conformally equivalent}, and 
for the classical series in finite dimension over $\K=\R$ a classification
has been given. Roughly, one considers
 the set of all $\alpha \in \Aut(\g)$
such that $F^* \circ \alpha$ is a grading-reversing involution
(where $F^*$ is conjugation by $F$); it is called the {\it structure
variety of $\Herm(A,*)$, resp.\ of $\Aherm(A,*)$} (cf.\ [Be00, Section
IV.2]). It contains, for
instance, all ``modifications'' or ``isotopes'' given by
$$
\alpha =\pmatrix{0& H \cr H^{-1} & 0 \cr},  \eqno (8.5)
$$
where $H$ is an invertible element in $\Herm(A,*)$.
Then one has to classify $G$-orbits in the structure variety.
In finite dimension over the reals, topological connected components of the
structure variety are homogeneous under $G$, and thus the task is
relatively easy. In infinite dimension, or over other base fields or
-rings, it seems possible that continuous families of non-isomorphic
 modifications may exist. This is an interesting topic for future 
research.

\msk \nin {\bf 8.7. Example: algebras of endomorphisms.}
Let $V$ be a $\K$-module equipped with a bilinear
symmetric or skew-symmetric form $b:V \times V \to \K$ which is non-degenerate 
in the sense that the map 
$\sigma \: V \to V^* := \Hom(V,\K), v \mapsto b(v,\cdot)$ is 
bijective. Let $A = \End(V)$ 
 and define for $X \in \End(V)$ the adjoint $X^* \in \End(V)$ by 
$X^*.v := \sigma^{-1}(\sigma(v) \circ X)$. 
Of course, in a topological context one has to add further assumptions in order
to ensure that $A$ is a c.i.a.\ and that the adjoint map is continuous;
e.g., one may assume that $\K$ is a topological field and $V$ finite-dimensional
over $\K$, or that we are in a Hilbert-space setting.
Then $\Phi_1(X)$ is the adjoint of $X \in \End(V \oplus V)$ w.r.t.
the bilinear form on $V \oplus V$ given by 
$$
\pmatrix{0 & -b \cr b & 0 \cr}.
\eqno (8.6)
$$
In particular, if $b$ is a scalar product over $\K=\R$, then
$\sp(A,*)$ really is the symplectic Lie algebra  $\sp(V \times V,\R)$. 
This is essentially the context considered in
[Be96] (see also [Be00, Ch. VIII.4]).
As is seen by elementary Linear Algebra (cf.\ loc.\ cit.), $\tilde \Phi_1:\P
\to \P$ is then the ``orthocomplement map'' with respect to (8.6) (where
$\P$ is the Grassmannian of subspaces of type $V$ in $V \oplus V$ having
complement of type $V$), and hence
the hermitian projective line corresponds to
the ``Lagrangian variety with respect
to the symplectic form'', and the anti-hermitian projective line corresponds to
``Lagrangians with respect to the quadratic neutral form''
into which the orthogonal group $\OO(V,b)$ can be imbedded.

\sectionheadline
{9. A quantum mechanical interpretation}

As explained in the introduction, there
is a strong structural analogy between the mathematics considered
in this work and the axiomatics of quantum mechanics.
In the following, we give some examples for this structural analogy
by proposing a ``dictionary"  between
the language of generalized projective geometries 
and the language of quantum mechanics.
This dictionary is by no means complete -- we do not attack 
 topics such as spectral theory of our observables or the
use of unbounded operators. However, it seems that the theory
of Jordan pairs and -triple systems is rich and flexible enough
to incorporate such aspects; we intend to investigate these questions
in future work.
Our  references for classical, linear Quantum Theory are
[Th81] and [Va85]. According to [Th81, p. 33], the ``Basic
Assumption of Quantum Theory'' is formulated as follows:
 ``The observables and states of a system are described
by hermitian elements $a$ of a $C^*$-algebra $A$ and by states on $A$.''
Let us see what this assumption implies if one tries to interprete
it on the level of the projective completion of the algebra of
hermitian elements. Consequently, 
we will start with the observables and not with the states.

\msk \nin
{\bf 9.1. Observables.}
The {\it space of observables} is the space $X^+$ of 
a generalized projective geometry $(X^+,X^-)$.
The space $X^-$ may be called the ``space of non-observables"
or the ``space of observers".
As {\it standard model} we may take the hermitian projective line
$X^+=\P_h$ over an (infinite-dimensional) associative involutive c.i.a
$(A,*)$. 
In this case, $X^+$ and $X^-$ are canonically isomorphic (the isomorphism
is a {\it canonical null-system} in the sense of [Be03a]).
For a general approach, it seems not necessary to assume that $\K=\C$.

\msk \nin
{\bf 9.2. States and pure states.}
A {\it state} is an {\it intrinsic subspace of $X^+$}, i.e.
a subset $Y \subset X^+$ which appears {\it linearly} (i.e. as an affine 
subspace)
with respect to {\it any} affinization $y \in X^-$.
Such subspaces correspond to {\it inner ideals of $V^+$} in Jordan theory
(cf.\ [Be02, 2.7.(4)], [BL04]). 
A {\it pure state} is an {\it intrinsic line}, i.e.\
 a proper intrinsic subspace
which is minimal for inclusion. The {\it superposition} of two
pure states is the intrinsic subspace generated by the two lines.
Under some additional assumptions,
pure states correspond to {\it division idempotents} of the Jordan pair,
and spaces of certain states form again a generalized projective geometry
(cf.\ [Ka01] for results that point into this direction).
Pure states correspond to 
{\it rank-one elements} (cf.\ [Lo94] for the notion of ``rank''),
 and they are closely related to
{\it chains} in the sense of Chain Geometry (cf.\ [H95]).

\msk \nin
{\bf 9.3. The Hamilton operator.}
A {\it Hamilton operator} is a polarity $p:X^+ \to X^-$
(cf.\ Section 6.1). A Hamilton operator is called {\it free}
if the polarity $p$ is an {\it inner polarity} in the sense of [Be03a]. 
In the standard model, there exists a free Hamilton operator $p_0$,
given by the matrix $F$ (called the ``natural polarity'' in
Section 8.6). Then a general  Hamilton operator
can be seen as a deformation or {\it modification} of the free
one as explained in Section 8.6;
in particular, via Equation (8.5) every
invertible hermitian element $H$ leads to new Hamilton operator
that needs not be conjugate to $p_0$. 
Note that the canonical identification $X^+ = X^-$ ($=\P_h$) in the standard
model is {\it not} a Hamilton operator because it is a null-system.

\msk \nin
{\bf 9.4. The time dependent Schr\"odinger equation.}
The {\it time dependent Schr\"odinger equation} is a dynamical differential
equation canonically associated to the Hamilton operator $p$.
Of course, here one thinks first of the {\it geodesic  differential equation}
in the symmetric space $M:=M^{(p)} \subset X^+$ associated to the
Hamilton operator $p$. (Every symmetric space carries a canonical
torsionfree connection $\nabla$ (cf.\ [Be03b]), and a geodesic is simply a smooth
map $\alpha: \K \supset I \to M$ which is compatible with connections.
In a chart, the geodesic equation is as usual $\alpha''(t)=C_{\alpha(t)}
(\alpha'(t),\alpha'(t))$ 
where $C$ is the Christoffel tensor of $\nabla$ in the chart.)
However, as pointed out in [AS97], the Schr\"odinger evolution should
rather be seen as a Hamiltonian flow and not as a solution of a second order 
differential equation. But it is possible to reconcile these two aspects
inside the category of generalized projective geometries because 
the tangent geometry $(TX^+,TX^-,Tp)$ is again of the same type,
and here the geodesic flow of $M^{(p)}$ appears as flow of a vector
field, namely of the {\it spray associated to the canonical connection
of $M^{(p)}$} (cf.\ [Be03b]). 

\msk \nin
{\bf 9.5. The time independent Schr\"odinger equation.}
An {\it eigenstate} of the Hamilton operator $p$ is an intrinsic
line which at the same time is a geodesic on $M^{(p)}$.
They correspond to division tripotents of the Jordan triple system
associated to $p$. A {\it complete system of eigenstates} corresponds
to a {\it frame} of the Jordan triple system. The {\it time
independent Schr\"odinger equation} consists in decomposing
a given tripotent with respect to a frame. 

\msk
\nin {\bf 9.6. Quantization.}
Note that some models of special and general relativity such as
Minkowski space and the de Sitter- and anti-de Sitter model
(and more general {\it causal symmetric spaces}) can be realized
via generalized projective geometries ([Be96], [Be00]).
It would be tempting to interprete a quantization of such spaces as
a sort of representation of these finite-dimensional geometries
in an infinite-dimensional geometry.

\sectionheadline{10. Prospects}

\nin {\bf 10.1. Generalizations.}
The differential calculus developed in [BGN03] works in more general
contexts, called ``$C^0$-concepts'', than the one of topological rings and
 modules. For instance, we may consider the class of rational mappings
defined on Zariski-open sets in finite-dimensional vector spaces over
an arbitrary infinite field $\K$ and define the class $C^1$ as in Section
1.3, where now $C^0$ means ``rational''. Essentially all results of the
present work carry over to this more general framework (details are left
to the reader).
In particular, all finite-dimensional Jordan algebras, -triple systems
and -pairs over arbitrary infinite fields are automatically 
``continuos (quasi-) inverse'' since the formulas for (Bergman-) inversion
clearly are rational. Thus, in finite dimensions over infinite fields,
the projective completion is always a ``smooth rational manifold'' in the
sense of [BGN03], and our construction yields ``smooth
rational symmetric spaces''. All notions of differential geometry from
[Be03b]  continue to make sense in this setting.

\msk \nin 
{\bf 10.2. Lie group actions.}
In the context of Theorem 5.3, one would like the projective group $G$
to be a Lie group acting smoothly on the projective completion $X^\pm$.
However, in general it seems impossible to define a Lie group structure
on $G$ because $G$ is defined by generators, and it is very hard to
find a good atlas for the subgroup $H$. In the real or complex Banach
set-up, this problem can be avoided by taking instead of $G$ and $H$
the ``much bigger'' groups $\Aut(\g)$ and $\Aut(\g,D)$ which are
Banach Lie groups, and then realizing $X^+$ as a quotient manifold under
the action of $\Aut(\g)_0$. This is the strategy used in [Up85]; it
needs a fair amount of non-trivial functional analysis and
does not carry over to more general situations. 

Nevertheless, the problem remains wether in our general set-up it is possible
to find some extension of $G$ to a Lie group $\tilde G$ acting smoothly
on $X^\pm$. For instance, in the case of the standard models (Sections
7 and 8) this is the case: in case of the projective line we may take
$\tilde G = \GL_2(A)$ which is indeed a Lie group (if $A$ is a c.i.a., then
the algebra $M_n(A)$ of $n\times n$-matrices with entries in $A$ 
is a c.i.a (cf.\ [Bos90], [Gl02]), and hence
$\Gl_n(A)$ is a Lie group), and in case of the (anti-) hermitian projective
line we may take $\tilde G = \Sp(A,*)$, resp.\ $\tilde G = \UU(A,A,*)$
which are unitary groups associated to an involutive c.i.a. and hence, as
we have seen in Section 8.6,  are Lie groups.
We intend to investigate the problem of Lie group extensions of
general projective groups in future work.

\def\entries{ 

\[AS97 Ashtekar, A. and T.A. Schilling, {\it Geometrical Formulation
of Quantum Mechanics}, preprint 1997, arXive:
gr-qc/9706069



\[Be96 Bertram, W., {\it On some causal and conformal groups}, J. Lie
Theory {\bf 6} (1996), 215 -- 247

\[Be00 Bertram, W., ``The Geometry of Jordan and Lie Structures,''
Lecture Notes in Math. {\bf 1754}, Springer-Verlag, Berlin, 2000


\[Be02 ---, {\it Generalized projective geometries: general
theory and equivalence with Jordan structures}, 
Advances in Geometry {\bf 2} (2002), 329--369  

\[Be03a ---, {\it The geometry of null-systems, Jordan algebras
and von Staudt's theorem}, Ann. Inst. Fourier {\bf 53} (1) (2003),
193 -- 225 

\[Be03b ---, {\it Differential geometry over general base fields
and -rings. Part I: First and scond order theory}, preprint, Nancy 2003

\[BGN03 Bertram, W., Gl\"ockner, H. and K.-H. Neeb, {\it 
Differentiable calculus, manifolds and Lie groups over arbitrary
infinite fields}, to appear in Expos. Math., arXive:
math.GM/0303300

\[BL04 Bertram, W. and H. L\"owe, {\it Incidence geometries associated
to Jordan pairs}, in preparation

\[BN03 Bertram, W., and K.-H. Neeb, {\it 
Projective completions of Jordan pairs. Part I. Geometries associated
to $3$-graded Lie algebras}, to appear in: J. of Algebra,
arXive: math.RA/0306272


\[BlHa01 Blunck, A. and H. Havlicek, {\it The connected components
of the projective line over a ring}, Advances in Geometry {\bf 1} (2001),
107--117


\[Bo90 Bost, J.-B., {\it Principe d'Oka, $K$-theorie et syst\`emes dynamiques 
non-commu\-ta\-tifs}, 
Invent. Math. {\bf 101} (1990), 261--333 

\[CGM03 Cirelli, R., Gatti M. and A. Mani\`a, {\it The pure state
space of quantum mechanics as hermitian symmetric space},
J. of Geometry and Physics {\bf 45} (2003), p. 267 -- 284

\[DNS89 Dorfmeiser, J., E. Neher, and J. Szmigielski, {\it
Automorphisms of Banach manifolds associated with the KP-equation}, 
Quart. J. Math. Oxford (2) {\bf 40} (1989), 161--195 

\[DNS90 ---, {\it Banach
manifolds and their automorphisms associated with groups of type
$C_\infty$ and $D_\infty$}, in ``Lie algebras and related topics,''
Proc.\ Res.\ Conf., Madison/WI (USA) 1988, Contemp.\ Math.\ {\bf 110} (1990),
43--65


\[FK94 Faraut, J. and A. Koranyi, {\it Analysis on Symmetric Cones},
Clarendon Press, Oxford 1994

\[Gl01a Gl\"ockner, H., {\it infinite-dimensional Lie groups without 
completeness 
restrictions}, in: Strasburger, A. et al. (Eds.)
``Geometry and analysis on Lie groups,'' Banach Center Publications,
Vol. {\bf 55}, Warsawa 2002; 53--59


\[Gl02 ---, {\it Algebras whose groups of units are Lie groups}, 
Studia Math. {\bf 153:2} (2002), 147--177 


\[dlH72 de la Harpe, P., ``Classical Banach Lie Algebras and
Banach-Lie Groups of Operators in Hilbert Space,'' Lecture Notes in
Math.\ {\bf 285}, Springer-Verlag, Berlin, 1972 

\[Ha82 Hamilton, R., {\it The inverse function theorem of Nash and
Moser}, Bull of the Amer. Math. Soc. {\bf 7} (1982), 65--222

\[H95 Herzer, A., {\it Chain geometries}. In: F. Buekenhout (editor),
``Handbook of Incidence Geometry,'' Elsevier 1995

\[IM02 Isidro, J. M. and M. Mackey, {\it The manifold of finite rank
projections in the algebra ${\cal L}(H)$ of bounded linear operators},
preprint, Santiago 2002

\[JNW34 Jordan, P., von Neumann, J. and E. Wigner, {\it On an algebraic
generalization of the quantum mechanical formalism}, Ann. Math. {\bf 35}
(1934), 29 -- 64



\[Ka83 Kaup, W., {\it \"Uber die Klassifikation der symmetrischen
her\-mi\-te\-schen Mannigfaltigkeiten unendlicher Dimension I, II},
Math. Annalen {\bf 257}(1981), 463--486; {\bf 262}(1983), 57--75

\[Ka01 ---, {\it On Grassmannians associated with $JB^*$-triples},
Math. Z. {\bf 236} \break (2001), 567--584 

\[Ke74 Keller, H. H., ``Differential Calculus in Locally Convex
Spaces,'' Lecture Notes in Math., Springer-Verlag, 1974 


\[Lo69 Loos, O., ``Symmetric Spaces I,'' Benjamin, New York, 1969

\[Lo75 ---, ``Jordan Pairs,'' Springer LNM 460, Berlin, 1975


\[Lo94 ---, {\it Decomposition of projective spaces defined by
unit-regular Jordan pairs}, Comm. Alg. {\bf 22} (10) (1994),
3925--3964


\[Lo96 ---, {\it On the set of invertible elements in Banach--Jordan algebras}, Results in Math. {\bf 29} (1996), 111-114 

\[MM01 Mackey, M., and P. Mellon, {\it Compact-like manifolds
associated to $JB^*$-triples}, manuscripta math. {\bf 106} (2001),
203--212 



\[MC03 McCrimmon, K., ``A Taste of Jordan Algebras,'' Springer-Verlag,
New York 2003


\[Ne02 Neeb, K.-H., {\it A Cartan-Hadamard theorem for Banach-Finsler
Manifolds}, Ge\-om. Dedicata {\bf 95} (2002), 115--156


\[PS86 Pressley, A., and G. Segal, ``Loop Groups,'' Oxford University Press, 
Oxford, 1986

\[Th81 Thirring, W., ``A course in Mathematical Physics. Vol. 3:
Quantum Mechanics of Atoms and Molecules,''
Springer-Verlag, New York 1981


\[Up85 Upmeier, H., ``Symmetric Banach Manifolds and Jordan
$C^*$-algebras,'' North Holland Mathematics Studies, 1985 

\[Va85 Varadarajan, V.S., ``Geometry of Quantum Theory,'' Springer-Verlag,
New York 1985

}

\references
\dlastpage 

\vfill\eject

\bye